\documentclass[12pt]{thloria}

\usepackage{cmll,aeguill}
\usepackage[T1]{fontenc}
\usepackage[latin1]{inputenc}
\usepackage{latexsym,multicol}
\usepackage{amsmath,graphics,stmaryrd}
\usepackage{amssymb}
\usepackage{amsfonts,verbatim}
\usepackage{graphicx,lscape}
\usepackage{epsfig,color}
\usepackage{array}
\usepackage{flafter}
\usepackage[all]{xy}

\newtheorem{theoreme}{Théorème}
\newtheorem{prop}{Proposition}

\newtheorem{lemma}{Lemme}
\newtheorem{definition}{Définition}
\newenvironment{proof}[1][Preuve]{\begin{trivlist}
\item[\hskip \labelsep {\bfseries #1}]}{\end{trivlist} \vspace{-0.5cm}
  \hspace{14cm} $\Box$}

\newcommand{\sep}{\ | \ }
\newcommand{\para}[1]{\noindent {\bf #1.}}

\newcommand{\Implies}{\Rightarrow}
\newcommand{\equi}{\Leftrightarrow}

\newcommand{\Par}{\parr}
\newcommand{\bool}{\mathbb{B}}

\newcommand{\prefixe}{\sqsubseteq}

\newcommand{\cat}{\mathbb{C}}
\newcommand{\path}{\twoheadrightarrow}
\newcommand{\p}{\mathcal{C}_G}
\newcommand{\n}{\mathcal{N}}
\newcommand{\A}{\mathcal{A}}
\newcommand{\B}{\mathcal{B}}
\newcommand{\C}{\mathcal{C}}
\newcommand{\D}{\mathcal{D}}

\newcommand{\chemin}{\twoheadrightarrow}

\newcommand{\paths}{\mathrm{Path}}

\newcommand{\Neg}{*}

\newcommand{\pop}{\multimap}
\newcommand{\tensor}{\otimes}
\newcommand{\soit}{\ \mathrm{new} \ }
\newcommand{\dans}{\ \mathrm{in} \ }
\newcommand{\si}{\ \mathrm{if} \ }
\newcommand{\alors}{\ \mathrm{then} \ }
\newcommand{\sinon}{\ \mathrm{else} \ }
\newcommand{\Unit}{\mathrm{Unit}}
\newcommand{\Bool}{\mathrm{Bool}}
\newcommand{\Nat}{\mathrm{Nat}}
\newcommand{\Ref}[1]{\mathrm{ref}[#1]}
\newcommand{\Skip}{\mathrm{skip}}
\newcommand{\zero}{\mathrm{zero}}

\newcommand {\sem}[1]{\ensuremath{\llbracket #1 \rrbracket}}

\newcommand{\rules}[3]{[#1] \ \frac{#2}{#3}}

\newcommand{\reduct}{\Downarrow}
\newcommand{\cone}{\mathfrak}
\newcommand{\point}{\mathnormal{2}}

\ThesisTitle{De l'opérateur de trace dans les jeux de Conway}
\ThesisAuthor{Nicolas Tabareau}

\DontShowLogos

\ThesisDate{2004-2005}
\ThesisChevaleret

\ContinuousNumbering{figure}
\NoChapterPrefix
\NoChapterNumberInRef

\begin{document}

\MakeThesisTitlePage

\tableofcontents

\DontWriteThisInToc

\chapter*{Résumé}

La compréhension sémantique des langages avec références est encore
assez partielle et l'on est loin d'une extension de l'isomorphisme de
Curry-Howard aux langages de programmation impératifs. Néanmoins, dans
un article récent Samson Abramsky, Kohei Honda et Guy McCusker
définissent une sémantique des jeux complètement adéquate (fully
abstract) pour un langage impératif en appel par valeur, avec
variables muables et types références~\cite{ahm}. Le modèle est défini au
moyen de jeux d'arènes de Hyland et Ong, et de stratégies filaires et
bien parenthésées.  Dans une note de recherche un peu plus
ancienne~\cite{milner:ac5}, Robin Milner interprète l'opérateur
``new'' du pi-calcul au moyen d'un opérateur de trace (ou de
``feedback''). Notre travail a pour but initial de comprendre ensemble
ces deux articles. Mais une autre motivation doit s'ajouter à ce
travail, motivation qui s'inscrit dans un vaste projet de refonte de
l'informatique et de sa théorie dans un cadre totalement algébrique.

Dans un premier temps, nous avons construit un modèle de
sémantique des jeux parenthésés de logique linéaire intuitionniste
disposant de plus d'un opérateur de trace. Pour cela, nous 
avons utilisé le modèle des jeux de Conway, augmenté avec une notion
de gain définie de manière axiomatique. Cette politique de définition
algébrique nous permet de rapprocher cette notion de gain d'une forme
de distance ou de norme sur l'espace des positions du jeu. Cela
participe d'une géométrisation du problème qui pourrait déboucher vers
un lien futur avec la Géométrie des Interactions.

Ensuite, nous avons montré l'existence d'un comonoïde commutatif
libre sur les jeux de Conway à gain. Cette façon
de voir permet de lever le voile qui masque bien souvent les
différentes définitions de l'exponentielle dans un même cadre
sémantique. Cette modalité acquiert ainsi un rôle plus intrinsèque et
ne fait plus l'objet de controverse.

Enfin, nous avons utilisé notre nouveau cadre pour décrire un modèle
d'un langage de type Algol avec fonctionnelle d'ordre supérieur. Pour
cela, nous avons allié le pouvoir de la logique linéaire pour
décrire l'aspect fonctionnel du langage ainsi que le pouvoir des
traces pour décrire l'aspect mémoriel du langage.

\mainmatter

\chapter{Carnet de bord} \label{ch:intro}

Dans cette introduction en forme de journal de voyage, nous allons
retracer les grandes lignes de notre réflexion sur la modélisation des
références dans un cadre algébrique. 

\section*{Prologue}

\para{Algèbre et sémantique}
Le projet global qui nous anime est de rattacher algèbre et
sémantique, le plus souvent grâce à des intermédiaires catégoriques. Ceci
permettra d'utiliser le formidable éventail de concept 
développé depuis plus de deux siècles dans ce champ mathématique. C'est d'ailleurs
un vaste projet qui anime nombre de mathématiciens, logiciens ou
informaticiens actuels comme Girard, Abramsky ou Hyland pour ne citer
qu'eux. Il est en effet 
inévitable que le séisme qui s'est produit dans la physique du début
du vingtième siècle se renouvelle en informatique sous l'impulsion
d'une vraie jonction avec l'algèbre. Si l'informatique peut être vue
comme la mathématique du discret, reste à donner un sens précis à
cette intuition. Déjà bon nombre de concepts catégoriques ont pris du
sens dans le cadre sémantique mais ce n'est pas suffisant pour
considérer qu'un pont solide a été établi. C'est par exemple pour le
renforcer que l'avant-garde de la communauté concurrente a introduit le
concept d'homotopie pour décrire les phénomènes d'interférences entre
deux processus s'exécutant en parallèle \cite{EGoubault:mscs2000}. Le
concept d'homotopie a 
aussi été le moyen de décrire des notions techniques comme l'innocence
en sémantique des jeux \cite{mellies:ag2}. Nous ne reviendrons pas ici
sur ces notions mais nous nous efforcerons de donner d'autres pistes
algébriques. 
Nous nous intéresserons au lien entre la notion catégorique d'un
opérateur de trace et le concept de référence dans les langages de
programmation; nous introduisons une notion positionnelle de 
gain sous la forme d'une distance entre deux positions; 
et décrivons différentes méthodes algébriques pour construire le monoïde
commutatif libre dans une catégorie.
\newline

\para{Trace et référence: lorsque $entr\acute{e}e = sortie$ fait sens}
Nous pensons que les traces qui formalise essentiellement le concept
sémantique de \emph{feedback} (lorsque l'entrée peut être ``branchée''
sur la sortie) forment un cadre intéressant pour décrire
les références ou les variables locales. En effet, une variable locale
peut être vue à la fois comme une entrée (l'écriture) et
une sortie (la lecture) qui vivent ensemble grâce à la trace.
Mais il nous faut donner vie à cette idée. La 
démarche adoptée est la suivante; regarder un modèle sémantique des
jeux de logique linéaire pour lequel on peut définir une trace. En
effet, l'aspect sémantique des jeux permet de décrire fidèlement la
partie fonctionnelle du langage (beaucoup de résultats récents de
``full abstraction'' passent par ce chemin) et l'existence d'une trace
permet d'ajouter la couche mémorielle nécessaire à la description des
références. 

La tâche semble malaisée tant la logique linéaire représente un moyen
de comprendre le séquentiel, ie. l'absence de boucle tandis que
l'existence de feedback est précisément un moyen d'en créer à loisir!
Pourtant, cela ne semble pas absurde d'essayer de faire vivre ensemble
ces deux notions car cela s'inscrit dans un processus de
rapprochement de la logique linéaire et de la théorie des noeuds
initié par des gens comme Paul-André Melliès via la notion d'homotopie.

Dans cette optique, nous avons dans un premier temps étudié les jeux
asynchrones \cite{mellies:ag4} pour nous rendre compte que
l'existence de traces était condamnée par la notion de gain donnée
originellement. En effet, ce gain n'étant pas autodual, il brise 
la symétrie au c\oe ur de l'opérateur de trace (en empêchant la
possibilité d'avoir une simplification à gauche).
Il nous a donc fallu faire un lourd travail pour isoler la notion de
gain de la distinction $\tensor$/$\Par$ en logique linéaire,
distinction à l'origine de cette asymétrie.

Méthodologiquement, nous nous sommes tournés vers un modèle de
sémantique des jeux où la trace était déjà présente, les jeux de
Conway. Ce modèle admet une trace de manière canonique car il est
compact fermé, ie. que le tenseur y est autodual. Malheureusement,
ce modèle a été un peu délaissé par la communauté sémantique, et des
notions cruciales comme le parenthésage y sont absentes. Nous avons
donc voulu les introduire sans pour autant briser la structure
compacte close. À terme, nous envisageons des jeux asynchrones avec un
opérateur de trace comme synthèse entre les jeux de Conway et les jeux
asynchrones. 
\newline  

\para{Distance et gain: une formalisation du contrôle}
La notion de gain a été pour nous une façon de recomprendre le
parenthésage et, par ce biais, un moyen d'augmenter notre catégorie de
Conway avec du contrôle. Le contrôle représente en particulier le
fait de pouvoir forcer un programme à regarder systématiquement un, deux
ou tous les arguments de la fonction qu'il calcule. 
Toujours dans cette optique algèbre et sémantique,
nous en donnons ici une version axiomatique qui permet de rapprocher le
gain au concept de distance entre deux positions. La distance
exprimant ici le nombre de questions Joueur et Opposant ouvertes et
non répondues (appelées \emph{questions pendantes}) entre ces deux
positions. 

Cette approche singulière nous a permis à la fois de préserver la
structure compacte fermée des jeux de Conway mais aussi de rapprocher
ces derniers d'une autre tentative d'articulation entre algèbre et
sémantique, à savoir la Géométrie des Interactions.

Plus précisément, l'axiomatique proposée permet de considérer que le
gain d'un chemin~$s$ calcule le nombre de questions posées par Joueur
et Opposant dans $s$. En particulier, un seul coup peut poser
plusieurs questions, ce qui n'arrive pas dans une approche par
justification avec pointeurs.

Il est important de noter que notre définition est axiomatique
et ainsi capture diverses notions de gain qui permettent de formaliser
en particulier divers comportements du jeu booléen. On peut dans cette
optique définir plusieurs jeux booléens donnant différentes
interprétations au jeu 
$$
\bool_1 \tensor \bool_2 \longrightarrow \bool_3
$$
Une dans laquelle le premier coup d'Opposant dans $\bool_3$ force à
interroger à la fois $\bool_1$ et $\bool_2$, ou alors l'un des deux,
ou bien même aucun.

On peut rapprocher cette définition (orientée distance) à la notion de
norme dans un espace de Hilbert et ainsi voir se dessiner un pont
entre jeux asynchrones tracés avec gain et Géométrie de
l'Interaction, même si ce lien reste malheureusement fantasmatique. 
\newline

\para{Comonoïde libre et exponentielle: le phénomène de duplication}

Voici en quelques lignes les motivations amenant au problème du calcul
du comonoïde libre.

\begin{itemize}
\item
  La logique linéaire est la voie qui nous ouvre à l'algèbre
\item
  Dans ce cadre, le comonoïde commutatif $1 \leftarrow {!A}
  \rightarrow {!A} \tensor {!A}$ décrit la copie
\item
  Un état est un objet copiable dans notre catégorie
\item
  Comment calculer $!A$ librement 
\end{itemize}

Nous avons ensuite défini un cadre dans lequel on peut calculer
l'exponentielle comme un comonoïde commutatif libre. Ici encore, cette
approche permet de donner un statut algébrique à une modalité de
logique linéaire. Nous donnons dans un premier temps un cadre agréable
dans lequel l'exponentielle se calcule simplement comme une extension
de Kan. Ensuite nous mentionnons un résultat de Dubuc pour construire
le comonoïde libre dans un cadre plus général où certaines propriétés de
commutation aux limites sont relâchées. Et enfin, nous étendons ce
résultat au cas du comonoïde commutatif libre, cas qui nous
intéresse ici.

Comme cette construction est
beaucoup moins lisse que celle par extension de Kan, nous espérons
dans un avenir proche pouvoir la réinterpréter elle aussi en terme
d'extensions de Kan. Cette optique ouvre d'ailleurs la voie au
développement d'une théorie monoïdale s'appuyant sur la théorie des
opérades de May mais ce champ est encore à explorer plus en détail.

Pour passer d'un modèle de MELL à un modèle de logique linéaire
intuitionniste, il nous manque la construction d'un produit
cartésien.   
Cette construction étant impossible pour des jeux Conway
généraux~\cite{mellies:ag3}, nous avons restreint la catégorie étudiée à celle
des jeux de Conway négatifs, où nous avons alors défini le produit
comme la simple union des deux jeux.

On est alors essentiellement en présence d'un modèle de logique
linéaire tracé et il nous faut maintenant trouver un langage pour
exprimer notre idée première qui est que les références se modélisent
avec des traces.
\newline

\para{Un modèle de langage avec trace}
Toutes ces investigations nous permettent de construire un langage de
type Algol avec des références (sans aliasing) à la fois globale et
locale. Les références y sont interprétées comme des variables
présentent à la fois en entrée et en sortie, et la localité est
obtenue en traçant sur cette entrée/sortie.
 
Ceci constitue l'aboutissement de la première étape dans notre
programme d'algébraïsation de la sémantique des langages de
programmation. Nous espérons même déboucher à une extension de
l'isomorphisme de Curry-Howard pour des langages évolués.
\newline

Revenons maintenant un peu plus en détail sur la réflexion qui nous a
fait aboutir à ce travail.

\section*{Sur la trace des références.}

\noindent {\bf De la Trace \dots}
Les catégories monoïdales tracées~\cite{joyal-street-verity} 
ont été introduites par Joyal, Street et Verity afin de fournir une
description uniforme de divers constructions mathématiques ayant un
comportement cyclique. Parmi les constructions les plus notables, nous
citerons la fermeture des tresses en théorie des n\oe uds et
l'opérateur de trace en algèbre linéaire.
Elles devinrent rapidement populaires dans la communauté informatique
comme un moyen élégant pour exprimer la notion de boucles dans un
cadre catégorique.
Elles ont été extrêmement fructueuses dans ce champ, que ce soit pour
formaliser la formule d'exécution de la Géométrie des Interactions
~\cite{abramsky:retracing,abramsky-haghverdi-scott}, pour analyser
l'opérateur de point fixe en théorie des domaines \cite{hasegawa02},
ou pour offrir un modèle catégorique en concurrence et plus récemment
en physique quantique~\cite{abramsky-cook}.
Formellement, une catégorie monoïdale tracée est une catégorie
monoïdale balancée (catégorie monoïdale avec tressage [braiding] et twist)
munie d'un \emph{opérateur de trace}
\begin{equation}\label{equation/to}
Tr_X: \qquad \frac{X \tensor A \quad \longrightarrow \quad X \tensor B}
{A \quad \longrightarrow \quad B}
\end{equation}
qui associe à chaque morphisme $f:X\tensor A\longrightarrow X\tensor B$
un morphisme $Tr_X(f):A\longrightarrow B$,
soumis a une série d'axiomes de cohérence rappelée en
section~\ref{section/tmc}. 

\begin{figure}[h!]
\begin{center}
\input{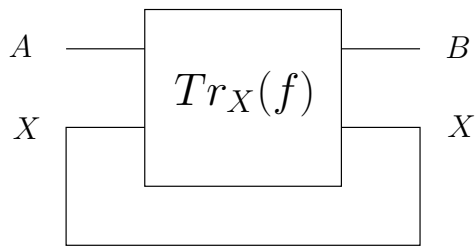}
\end{center}
\caption{Diagramme représentant l'action de la trace comme une
  redirection de la sortie vers l'entrée (essence du feedback)} 
\end{figure}

\para{\ldots aux références}

L'une des premières apparitions de l'opérateur de trace pour
modéliser le feedback dans un cadre sémantique est due à
Milner~\cite{milner:ac5}. Dans ce papier, il présente une façon
abstraite de modéliser le feedback dans les \emph{action calculi} qu'il
nomme \emph{réflexion} (le terme trace n'était pas encore à l'ordre du jour).
La réflexion lui permet de décrire des opérations compliquées comme
le célèbre opérateur de restriction $\nu :
\epsilon \rightarrow p$ qui devient simplement la trace de la
diagonale $(x)\langle xx\rangle: p \rightarrow p \tensor p$.
Rappelons que cette opérateur modélise le fait qu'un canal de communication public
peut être soudain restreint pour devenir un canal de communication
privée entre les processus qui communiquaient déjà dessus. On peut
recomprendre ce mécanisme comme le passage d'une mémoire globale à
une mémoire locale via un phénomène de localisation. Il est à noter
que bien que Milner n'avait pas connaissance des travaux récents sur
les traces à l'époque où il a défini les réflexions, tous
les axiomes qu'il donne pour que son opérateur de restriction conserve
les propriétés habituelles (comme le fameux ``scope extrusion'')
coïncident exactement avec l'axiomatique de l'opérateur de trace. Ceci
fait de la trace un objet canonique qui semble destiné à interpréter
les références.
  
Ici, nous nous intéressons à la trace comme moyen de description des
variables locales dans les langages de
programmation. Traditionnellement en sémantique, on 
interprète un langage de programmation dans une catégorie en
distinguant les objets $A$ décrivant les valeurs, des objets $TA$
décrivant les calculs de type $A$, où $T$ est une monade. Dans
le cas des références, la monade considérée est la monade d'état $S
\pop ( S \tensor \_)$ qui permet d'interpréter un programme de type $A\rightarrow
B$ comme un programme prenant une valeur $A$ et renvoyant un calcul
$S \pop ( S \tensor B)$, ce qui, via la clôture monoïdale, correspond à
un morphisme de $S\tensor A \rightarrow S\tensor B$. 

Il faut penser cette interprétation comme la description d'un
système avec entrée/sortie et mémoire accessible à l'utilisateur vu
comme un morphisme $f:S\tensor A \longrightarrow S\tensor B$.
Dès lors, si on est capable de prendre la trace sur $S$ de
$f$, on obtient la description avec mémoire interne à savoir les
morphismes usuels de type $A\longrightarrow B$. C'est l'analogue de la
restriction (ou localisation) chez Milner pour un langage où les
canaux sont remplacés par des adresses mémoires.

Par la suite, cette approche va être utilisée pour décrire un modèle
d'un langage de type Algol avec fonctionnelle d'ordre supérieur. Pour
cela, il nous faut allier le pouvoir de la logique linéaire pour
décrire l'aspect fonctionnel du langage ainsi que le pouvoir des
traces pour décrire l'aspect mémoriel du langage. 

\section*{Les jeux asynchrones ou la face nord de l'Éverest.}

Nous sommes partis des jeux
asynchrones~\cite{mellies:ag4} car c'est un modèle de logique linéaire où
des concepts de théorie des n\oe uds sont déjà intégrés
via la notion de chemins homotopes. 

Malheureusement, l'espoir d'y trouver une trace de manière directe a
été vite vain car il y a quelques problèmes difficiles à
surmonter. Ceci découle en particulier de la remarque suivante.
\newline

\para{Catégorie ponctuée}
Supposons qu'il existe un objet initial $0$ et un objet terminal
$\top$ au sein d'une catégorie symétrique monoïdale fermée.
Le foncteur $$A\mapsto A\tensor B$$ a un adjoint à droite et préserve
donc les colimites pour tout objet $B$ de la catégorie.

On en déduit qu'il existe un unique isomorphisme:

$$0\tensor B \cong 0$$

et plus généralement un unique isomorphisme 
$$0\tensor A \cong 0\tensor B$$
pour tout objet $A$ et $B$ de la catégorie.

En particulier, en instanciant avec $A=\top$ et $B=0$,
il y a un unique isomorphisme:

$$f:0\tensor \top \longrightarrow 0\tensor 0.$$

Maintenant, supposons que la catégorie est tracée.
On peut calculer la trace sur $0$ de~$f$

$$Tr_{0}(f):\top \longrightarrow 0.$$

Il suit la coïncidence de l'objet initial et terminal dans la
catégorie, modulo un unique isomorphisme.
Une telle catégorie est souvent appelé \emph{ponctuée}.
\newline

\para{Vers des jeux asynchrones tracés}
La catégorie des jeux asynchrones formulée dans \cite{mellies:ag3}
n'est \emph{pas} ponctuée car il n'existe pas de stratégie du jeu 
$\top$ dans le jeu $0$.
On en déduit qu'elle n'est pas tracée.

On va donc sortir pour le moment du cadre des jeux asynchrones pour
construire étape par étape un modèle possédant toutes les propriétés
annoncées plus haut. L'idée est ensuite de pouvoir revenir au cadre
des jeux asynchrones mais comme le souci d'algébraïsation des
outils utilisés dans la construction a amené un lourd travail, nous ne
pouvons présenter ce cadre ici.

\section*{La piste compacte close.}

L'origine de la notion de trace vient de son existence
automatique pour les catégories compactes closes. 
Si l'on voit la notion de trace comme la généralisation d'un monoïde
simplifiable (à gauche), la notion de catégorie compacte close est
alors la généralisation d'un groupe.

Pour illustrer ce concept dans l'univers mathématique, citons comme
exemple frappant la catégorie des espaces 
vectoriels avec le produit tensoriel usuel. Dans cette catégorie
symétrique monoïdale close, l'opérateur de clôture possède la
propriété remarquable d'avoir lui aussi une structure tensorielle. 
Un autre exemple naturel nous est donné par les catégories linéaires
dans lesquelles $\tensor = \Par$. Dans ces catégories, on sait
directement que le tenseur est autodual car
$$
(A \tensor B) ^\Neg = A^\Neg \Par B^\Neg = A^\Neg \tensor B^\Neg
$$

Plus généralement, les catégories compactes fermées sont symétriques
monoïdales closes avec un tenseur auto-dual. En d'autres termes, la
clôture est donnée par 
$$A \pop B \equiv A^\Neg \tensor B$$
Ainsi, tout morphisme $f:(X \tensor
A)^\Neg \tensor X \tensor B$  peut être transformé en
$\widehat{f}:(X^\Neg \tensor X)^\Neg \tensor A^\Neg \tensor B$ par
commutativité et associativité du tenseur, et par composition avec
l'identité, on obtient
$$Tr_X(f) = \widehat{f}(id_X)$$
On voit donc que toute catégorie compacte close est tracée. 
Mais n'oublions pas que nous voulons utiliser l'opérateur de trace
pour modéliser les références et l'on doit se demander si se
restreindre n'est pas trop fort au sens où l'on oublierait au passage
certaines catégories primordiales pour la description des références.

Un première réponse à cette question est la construction $Int$ due à
~\cite{joyal-street-verity} de la catégorie compacte libre engendrée
par un catégorie tracée. Cela dit que l'on peut toujours voir une
catégorie tracée comme une version laxiste d'une catégorie compacte
fermée. 
Cette construction est à la base de la notion de polarité en
logique.  

\begin{figure}[h]
\begin{center}
  \input{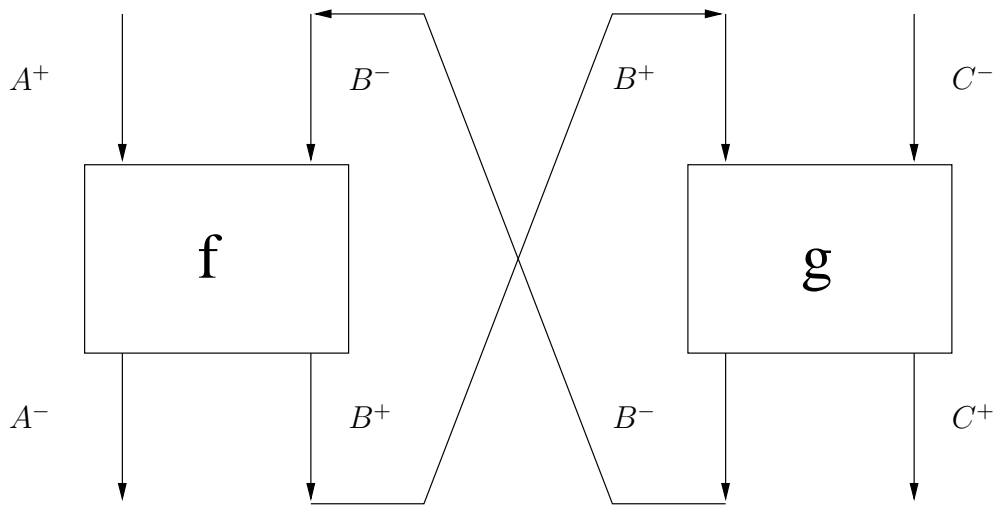}
\end{center}
\caption{Composition dans la catégorie $Int(\C)$}\label{fig:compositionInt}
\end{figure}

\begin{definition}
  Soit $\C$ une catégorie tracée. On définit $Int(\C)$ comme la
  catégorie compacte fermée ayant pour objet les couples $(A^+,A^-)$
  d'objets de $\C$ et pour morphismes entre $(A^+,A^-)$ et $(B^+,B^-)$
  les flèches $A^+ \tensor B^-
  \longrightarrow A^- \tensor B^+$ dans $\C$.  
  La composition est définie à l'aide de la trace (cf diagramme
  \ref{fig:compositionInt}). 
  Le dual est juste l'inversion de la polarité $(A^+,A^-)^\Neg =
  (A^-,A^+)$et le produit tensoriel est défini point à point
  $$
  (A^+,A^-) \tensor (B^+,B^-) = (A^+ \tensor B^+,A^- \tensor B^-)
  $$
\end{definition}

Du point des groupes et des monoïdes simplifiables, cette construction
correspond simplement au groupe libre sur un monoïde, l'objet $(x,1)$
représentant $x$ et l'objet $(1,x)$ représentant son inverse $x^{-1}$. 

Remarquons que cette construction est au c\oe ur de la catégorisation
des modèles de GoI donnée par Abramsky \cite{abramsky-haghverdi-scott}. 
La composition dans la catégorie $Int(\C)$ représentant la formule
d'exécution. 

Nous voulons pour nos travaux que le tenseur possède un adjoint à
droite donnant lieu à une clôture monoïdale, clôture qui n'est pas
reliée à priori à la clôture  
automatique des catégories compactes fermées. C'est en cherchant dans
cette voie que nous avons été amenés à regarder comment une catégorie
peut hériter de la clôture d'une catégorie qui la contient à
isomorphisme près.

\begin{prop}\label{prop:cloture}
Soit $(\C,\tensor_\C,\pop_\C)$ une catégorie symétrique monoïdale fermée
et $(\D,\tensor_\D)$ une catégorie symétrique monoïdale.
Supposons qu'il existe un foncteur monoïdal fort $U:\D \rightarrow \C$
qui est à la fois plein et fidèle et qui possède un adjoint à droite
$F:\C \rightarrow \D$ formant l'adjonction $U \dashv F$.

Alors, on peut exporter la fermeture sur $\C$ en une fermeture sur
$\D$ en définissant pour $A,B$ dans $\D$:
$$
A \pop_\D B = F(U(A) \pop_\C U(B))
$$
\end{prop}

Ce résultat stipule donc qu'une catégorie tracée $\C$ étroitement liée
à la catégorie $Int(\C)$ (via une adjonction) possède automatiquement
une structure close. C'est un premier résultat en faveur de l'étude
des catégories compactes closes, mais il manque une sorte de
réciproque à cette propriété, qui assurerait que si une catégorie
tracée $\C$ est close, alors elle est étroitement liée à $Int(\C)$. De
façon réjouissante, nos discussions avec Masahito Hasegawa (Juillet
2005) sur les relations qu'entretiennent la catégorie des jeux de
Conway négatifs $\n$ et sa ``compactifiée'' $Int(\n)$ l'ont amené à
préciser la situation dans le cadre des catégories tracées.

\begin{prop}
  Soit $\C$ une catégorie tracée. On peut définir un foncteur monoïdal
  fort plein et fidèle $J:\C \rightarrow Int(\C)$ par : $A \mapsto
  (A,I)$.On a alors:

  $\C$ est close ssi J admet un adjoint à droite
\end{prop}

Il est donc clair que si une catégorie tracée est close, ça clôture vient
directement de la clôture de $Int(\C)$, ie. de la catégorie compacte
fermée sous-jacente. Dès lors, il n'y plus de raison de s'empêcher de
travailler avec une catégorie compacte close, quitte à se restreindre
légèrement à posteriori.

\section*{La route du bagnard (Conway Game).}

Il semble à présent naturel de se pencher sur l'une des seules catégories
de sémantique des jeux ayant une structure compacte close: 
\emph{les jeux de Conway}. Cette catégorie a de plus le bon goût de
ressembler d'assez près aux jeux asynchrones et on peut donc espérer
rejouer le même scénario pour trouver une catégorie $*$-autonome.

En effet, les jeux de Conway ne sont pas beaucoup plus que des graphes
pour lesquels les idées d'\emph{innocence via classes d'homotopie}
et de structure de \emph{gain} pour permettre de retrouver une notion
de bon parenthésage semble pouvoir marcher. 
\newline 

\para{Remarque}
Nous avons étudié le cas des jeux de Conway asynchrones, en ajoutant
une notion d'homotopie sur les chemins pour donner une définition
algébrique de l'innocence. Nous avons aussi regardé la notion de
stratégie positionnelle (ie. qui se décrit par une relation sur les
positions) en montrant que toute stratégie innocente était
positionnelle. Dans ce cadre, nous voulons retrouver la trace des
stratégie positionnelle par la trace définie sur la catégorie des
relations. Malheureusement, il faut savoir sélectionner ce dont on
parle et donc ces résultats seront exprimés dans de futurs papiers. 
\newline

Étudions plutôt ici la notion de gain, pilier de l'algébrisation  du
contrôle,  car celle-ci peut facilement
briser la structure tracée, comme le montre la définition du gain pour
les jeux asynchrones. 

\section*{Le péage à double sens.}
Comme nous l'avons vu plus haut, la notion de catégorie monoïdale
tracée catégorise (et ainsi généralise) la notion habituelle de
monoïde simplifiable à gauche~$(M,\cdot,e)$.
De ce point de vue, l'opérateur de trace remplace l'implication:
\begin{equation}
\label{equation/lc}
\forall (x,a,b)\in M\times M \times M,
\hspace{2em}
x\cdot a = x \cdot b \ \Rightarrow \ a = b.
\end{equation}
Afin d'interpréter la logique linéaire propositionnelle et de
généraliser la condition de bon parenthésage, Paul-André Melliès a
assigné un gain $\kappa_A(x)\in\mathbb{Z}$ à chaque position $x$ d'un
jeu tout en demandant à ce que toute stratégie $\sigma$ ne joue que des
positions $x$ avec un gain positif: $\kappa_A(x)\geq 0$.
Malheureusement, le gain défini dans~\cite{mellies:ag3}
ne satisfait pas de propriété comme~(\ref{equation/lc}).
Plus précisément, étant donné trois jeux $X$, $A$ et $B$ et trois
positions $x$ de $X$ , $a$ de $A$ et $b$ de $B$, le gain $\kappa$ ne
vérifie pas:
$$
\Downarrow \qquad
\frac{\kappa_{X \tensor A \pop X \tensor B}
(x \tensor a \pop x \tensor b) \geq 0}
{\kappa_{A \pop B}(a \pop b) \geq 0}
$$
où $x\tensor a \pop x\tensor b$ et
$a \pop b$ représente, comme on si attend, les positions dans les jeux $X
\tensor A \pop X \tensor B$ et $A \pop B$ respectivement.
Ceci est embêtant car toute définition raisonnable d'un opérateur de
trace $Tr$ doit demander que la stratégie 
$$Tr_{X}(\sigma)  \ : \ A \longrightarrow B$$
joue la position
$$a\pop b$$
à chaque fois que
$$\sigma \ : \ X \tensor A \longrightarrow X \tensor B$$
joue la position
$$x \tensor a \pop x \tensor b$$
à partir d'une position $x$ du jeu $X$
sur lequel la trace est calculée.
Voilà pourquoi nous devons retoucher la notion de gain pour ne pas
briser la trace existant dans les jeux de Conway. En un sens, le gain
doit être complètement symétrique en Joueur et Opposant pour préserver
la structure compacte fermée. 
Nous définissons le gain $\kappa_A$ sur les chemins $s$ par une paire
d'entiers naturels $\kappa_A(s)\in\mathbb{N}\times\mathbb{N}$
représentant intuitivement le nombre de questions posées par Joueur
(première composante) et par Opposant (deuxième composante).

La condition habituelle de bon parenthésage est alors reformulée (et
généralisée) en demandant à ce que chaque chemin joué par une
stratégie satisfasse 
$$
\kappa_A^+(s) = 0 \Implies \kappa_A^-(s) = 0
$$
c'est-à-dire que si le Joueur est interrogé, il doit réagir en
répondant on en interrogeant à son tour.

\section*{L'ascension de la tour Exponentielle.}

Fort de ces premiers pas, il nous faut maintenant définir les
opérateurs de logique linéaire qui manque à notre modèle intuition.
Regardons d'abord l'exponentielle et essayons de voir quels liens elle
entretient avec les structures algébriques. 

A l'examen d'exemples tels que les espaces de cohérence ou les
catégories linéaires de Lafont, il semble que la bonne façon de
relier l'exponentielle à l'algèbre est de donner une construction
comonoïde commutatif libre. Nous nous sommes donc intéressés à un
objet d'étude courant en catégorie, la construction du monoïde libre.

Une construction du monoïde commutatif libre existe lorsqu'on a
la chance que le tenseur commute avec la somme. On doit simplement
calculer la limite
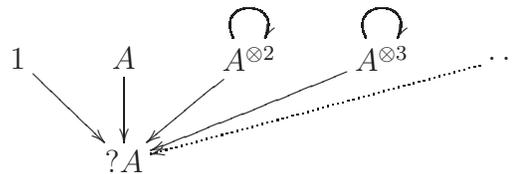
\vspace{0.2cm}
$$
\xymatrix{
1 \ar[rd] & A \ar[d] & A^{\tensor 2} \ar@(ul,ur)[] \ar[ld] &
A^{\tensor 3} \ar@(ul,ur)[] \ar[lld]&
\ldots \ar@{.>}[llld]\\
& {?A}  &&& }
$$

Malheureusement, ce n'est pas le cas dans les jeux de Conway car la
somme n'existe même pas en général et il faut trouver mieux.
Nous avons d'abord voulu recomprendre la construction susmentionnée
avec des extensions de Kan dans le but de pouvoir ensuite
généraliser.
Il s'avère qu'on obtient ce monoïde en calculant l'extension de Kan
sur $\mathcal{S}_0$ (la catégorie des ensembles fines et fonctions
ensemblistes) du jeu $A$ vu comme un foncteur monoïdal de $Bij$ (la
catégorie des ensembles fines et fonctions bijectives) dans la
catégorie $\C$ des jeux de Conway.  

$$
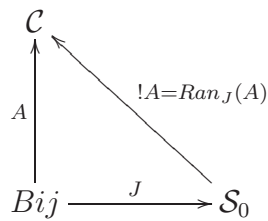
\xymatrix{
\C & & \\
& & \\
Bij \ar[uu]^{A} \ar[rr]^{J} & & \mathcal{S}_0 \ar[lluu]_{!A = Ran_J(A)}
}
$$

Pour généraliser ce premier résultat, nous nous sommes tournés vers une
construction du monoïde libre due à Dubuc qui dit essentiellement si
le tenseur commute uniquement aux colimites filtrées, on doit
calculer la même sorte de limite mais au lieu de prendre la catégorie
discrète, il faut prendre la catégorie filtrée simpliciale étendue aux
ordinaux.

Ceci convient très bien à notre catégorie de Conway car elle est
compacte close et le tenseur commute donc à toutes les limites et
colimites existantes (il possède un adjoint à droite et à
gauche). Mais nous voulons un comonoïde commutatif. Nous devons donc
étendre le résultat de Dubuc à la construction du monoïde commutatif
libre. 

Pour cela, nous nous sommes restreints aux catégories dont le tenseur
commute aux limites $\omega$-filtrées (cette restriction vient d'un
souci de simplicité mais la même idée semble s'appliquer à n'importe
quel type de limites $\varphi$-filtrées) et nous avons revisité la
construction de Dubuc avec des outils plus simples et modernes. Nous
sommes passés de la catégorie simpliciale à la catégorie des
injections pour récupérer les permutations. Ensuite, nous avons montré
que la colimite sur cette nouvelle catégorie était la monoïde
commutatif libre via une factorisation par le résultat de Dubuc.  

Ceci nous permet de présenter la construction du monoïde et du
comonoïde commutatif libre dans la catégorie des jeux de Conway à
gain via un cadre entièrement algébrique.

Nous avons ensuite voulu recomprendre ce résultat assez technique en
terme d'extensions de Kan mais cette partie de notre travail n'est pas
encore arrivée à son terme.
Cette dernière option semble néanmoins ouvrir la voie au développement d'une
théorie monoïdale à l'instar des théories algébriques de
Lawvere. Peut-être faut-il regarder du coté des opérades pour résoudre
ce joli problème.

Toujours est-il que nous avons pu utiliser la construction de Dubuc
pour construire notre exponentielle et ainsi obtenir un modèle de MELL
avec un forte assise algébrique.

\section*{Accommoding the additives}

Reste maintenant à définir les additifs. Ces
derniers ne peuvent pas s'obtenir sur la catégorie entière des jeux de
Conway comme il a été montré dans ~\cite{mellies:ag3}. 
Il faut alors se restreindre à la catégorie des jeux négatifs pour
pouvoir définir sereinement un produit cartésien, produit qui se
définit comme la juxtaposition des deux jeux. 
Le problème qui se pose alors est l'existence de l'adjoint à droite du
produit tensoriel. En effet, on ne peut plus utiliser la construction
du dual (qui donnerait un jeu positif) et on perd donc la structure
compact close. Heureusement, il existe un moyen d'exporter
automatiquement la fermeture via la proposition \ref{prop:cloture}. 

Il est à noter que même si cette propriété semble être folklorique en
catégorie (du moins était-elle connue par Martin Hyland), nous n'avons
pas pu trouver de références pour le moment. Cela semble indiquer que
ce fait remarquable reste méconnu dans le milieu sémantique.

Ici, l'adjonction servant de base à l'extension de la clôture est
entre le foncteur d'inclusion et le foncteur $Neg$
qui prend un jeu quelconque et oublie les parties commençant par des
coups joueurs.
 
Ainsi, il n'y a pas d'angoisse et on peut tout à fait travailler avec
la sous-catégorie des jeux négatifs tout en conservant la fermeture et
l'opérateur de trace.

\section*{Vers le Graal mémoriel}
Toutes ces considérations nous ont amené à l'élaboration d'un langage de
type Algol qui comporte de la mémoire à la fois locale et
globale. L'intérêt majeur de l'interprétation que nous en donnons réside
en deux points:
\begin{itemize}
\item
  \emph{sa simplicité.} Le cadre développé permet d'interpréter le langage
  sans ajouts byzantins venant rendre peu intuitif le résultat obtenu
\item
  \emph{son caractère algébrique.} La construction du modèle repose
  essentiellement sur des considérations catégoriques et donne ainsi
  à notre interprétation une portée plus générale. En effet, le schéma
  introduit peut s'appliquer à d'autres types de sémantique qui peuvent
  donner lieu à des interprétations inattendues et venant renforcer
  notre compréhension des langages de programmation. 
\end{itemize}

\section*{Travaux à suivre}

Dans les mois qui viennent, nous étudierons en autres les points
suivants:

\begin{itemize}

\item
{\bf Jeux de Conway asynchrones.}
Comme mentionné ci-dessus, nous voulons enrichir notre catégorie pour
retrouver le cadre des jeux asynchrones. Cela participe d'une
recompréhension des divers travaux en sémantique des jeux dans le
cadre unifié des jeux positionnels.

\item
{\bf Référence avec aliasing.}
Nous présentons ici un langage avec référence où les variables de
type $ref(ref(-)$ sont interdites. Cela est commode et suffisant pour
avoir un grand pourvoir d'expression mais nous ne pouvons
malheureusement pas parler de \emph{tas}. Ainsi, pour décrire cette
structure usuelle en informatique, il nous faudra lors de prochains
travaux relâcher cette contrainte d'anti-aliasing.
  
\item
{\bf Langages de bas niveau.}
Une partie de notre projet est de pouvoir donner une sémantique
uniforme lors du processus de compilation décrivant aussi bien le
langage haut niveau de l'utilisateur que l'assembleur ou le langage
machine en bout de chaîne. Ainsi, la sémantique pourra dépasser le
cadre usuel de test d'équivalence de deux programmes écrits dans le
même langage, et donnera par cette occasion le moyen de vérifier
jusqu'au bout la correction d'un compilateur réel.

\end{itemize}

\chapter{Jeu de Conway à Gain}\label{sct:payoff conway games}

Comme nous l'avons annoncé dans le chapitre précédent, nous allons
maintenant construire une catégorie de jeux compacte fermée avec une
notion de parenthésage intégrée.

\section*{Jeux de Conway.}

Avant de présenter notre modèle à gain, il parait nécessaire de
rappeler la définition des jeux de Conway tant ce formalisme a été
boudé par les sémanticiens des jeux. Il reste pourtant un des modèles
de jeux les plus naturels, et permet de rapprocher la sémantique de
notions algorithmiques en utilisant explicitement la structure de
graphes.   

Un jeu de Conway~\cite{joyal} $A = (V_A,E_A,\lambda_A)$ est la donnée :
\begin{itemize}
\item
  d'un graphe orienté enraciné  $(V_A,E_A)$ de racine $\star_A$ 
\item
  d'une fonction $\lambda_A: E_A \rightarrow \{-1,+1\}$ donnant la polarité d'un coup.
\end{itemize}
Comme d'habitude $-1$ signifie opposant et $+1$ joueur.

Il faut maintenant mentionner le vocabulaire usuel concernant les
parties d'un jeu.
\newline

\para{Chemins}
Comme d'habitude en théorie des jeux, on note $x \rightarrow y$
lorsque $(x,y) \in E_A$ et appelle chemin tout suite de coups $x_0
\xrightarrow{m_1} x_1 \xrightarrow{m2} \ldots \xrightarrow{m_{k-1}}
x_{k-1} \xrightarrow{m_k} x_k$. 
Dans ce cas, on note $m_1 \ldots m_k : x_0
\path x_k$ pour indiquer la position initiale et finale du chemin.
\newline

\para{Coup initial}
On appelle \emph{coup initial} d'un jeu de Conway $A$ toute flèche de
$E_A$ partant de la racine $\star_A$.  
\newline

\para{Parties}
Une \emph{partie} est un chemin partant de la racine $\star_A$
$$
\star_A \xrightarrow{m_1} x_1 \xrightarrow{m2} \ldots
\xrightarrow{m_{k-1}} x_{k-1} \xrightarrow{m_k} x_k 
$$ 
L'ensemble des parties est noté $P_A$.
\newline

\para{Parties alternées}
Une partie $m_1 \ldots m_k :\star_A \path x$ est dite \emph{alternée} lorsque:
$$
\forall i \in \{1 , \ldots , k-1\} \quad \lambda_A(m_{i+1}) = - \lambda_A(m_i)
$$ 

Nous avons défini les objets de la catégorie des jeux de Conway mais il
reste maintenant à exprimer les morphismes entre de tels objets. Ceci
est réalisé en ajoutant une notion de stratégie sur un jeu $A$ et en
donnant un peu de structure pour décrire ce qu'est une stratégie de
$A$ vers $B$ (et ceci de manière auto-duale pour avoir la structure
compacte close). 
\newline

\para{Stratégies}
Une stratégie $\sigma$ est un ensemble de parties alternées de longueur
paire tel que :  

\begin{itemize}
\item
  la stratégie $\sigma$ contient la partie vide $\epsilon$
\item
  toute partie non-vide commence par opposant
\item
  $\sigma$ est close par préfixe paire 
  $$
  \forall s \in P_A \ \forall m,n \in E_A \quad s \cdot m \cdot n \in \sigma
  \Implies s \in \sigma
  $$
\item
  $\sigma$ est déterministe : $\forall s \in P_A \ \forall m,n,n' \in
  E_A$
  $$
  s \cdot m \cdot n \in \sigma \mbox{ et } s \cdot m \cdot n' \in
  \sigma \ \Implies \ n=n' 
  $$
\end{itemize}

On note $\sigma : A$ lorsque $\sigma$ est une stratégie de $A$.

La plus petite stratégie est la stratégie $\{\epsilon\}$, qui ne
répond jamais.
On l'appelle la \emph{stratégie vide}, notée $\bot$.
\newline

\para{Un peu de structure}
Le dual d'un jeu de Conway $A$ est le jeu 
$$A^\Neg = (V_A , E_A ,
-\lambda_A)$$

Le produit tensoriel de deux jeux $A$ et $B$, noté $A
\tensor B$ : 

\begin{itemize}
\item
  $V_{A \tensor B} = V_A \times V_B$ 
\item
  $
  x \tensor y \rightarrow 
  \left\{
  \begin{array}{l}
    x' \tensor y \mbox{ if } (x,x') \in E_A \\
    x \tensor y' \mbox{ if } (y,y') \in E_B \\
    \end{array}
  \right.
  $  
\item
  $
  \begin{array}{l}
    \lambda_{A \tensor B}( (x \tensor y) \rightarrow (x' \tensor y)) = \lambda_A(
    x\rightarrow x')\\ 
    \lambda_{A \tensor B}( (x \tensor y) \rightarrow (x \tensor y')) = \lambda_B(
    y\rightarrow y')  
  \end{array}
  $
\end{itemize}

Il est notable que le tenseur réalise essentiellement le produit des
deux graphes sous-jacents. Ainsi, toute partie du produit tensoriel
peut être vue comme l'entrelacement de deux parties dans $A$ et dans
$B$. 
Le jeu de Conway $1 = (\{ \star \},\emptyset,\lambda$ est évidemment
l'élément neutre de cette loi monoïdale.

Du point de vue logique linéaire, on devrait définir la loi monoïdale
$\Par$, duale du tenseur, afin d'avoir l'implication linéaire $\pop$
et de définir les morphismes entre deux jeux. Mais comme nous voulons
une catégorie compacte close, le tenseur est ici auto-dual, ce qui
implique que $\tensor = \Par$.

On veut maintenant définir un morphisme de $A$ vers $B$ comme un
stratégie de $A^\Neg \tensor B$. Pour cela, il faut avoir une notion
d'identité et de composition de deux telles stratégies.
\newline

\para{Identité}
Classiquement, on considère la stratégie d'imitation (\emph{copycat
strategy}) de type $A^\Neg \tensor A$ comme la stratégie identité de
$A$ dans $A$. Dans ce qui suit, on note $A_1$ et $A_2$ pour
distinguer entre les deux copies de $A$.
$$
id_A = \{ s \in P_{A \pop A} \ | \ \forall s' \prefixe^{even} s
\quad  s'_{|A_1} = s'_{|A_2} \} 
$$

\para{Interactions}
On dit que $u$ est une interaction de $A,B,C$, notée $u \in int(A,B,C)$
si la projection de $u$ sur chaque jeu $A ^\Neg \tensor B$,$B ^\Neg
\tensor C$ et $A ^\Neg \tensor C$ est une partie.

Nous mettons maintenant en place le cadre pour notre définition de
stratégie gagnante. La condition de gain pour une stratégie est locale
car elle porte non seulement sur les parties jouées mais encore sur
les chemins apparaissant dans l'interaction. C'est ainsi que la notion
de contrôle s'initie dans notre cadre car un programme peut maintenant
être contraint à chacun de ses coups et non plus uniquement dans une
interaction globale. 
\newline 

\para{Chemin joué par une stratégie}
Un stratégie $\sigma$ \emph{joue un chemin} $t:x \path y$
lorsqu'il existe une partie $s:\star \path x$ dans $\sigma$ telle que
la composition $s;t:\star \path y$ est aussi dans $\sigma$. On
dit aussi dans ce cas que le chemin $t$ est dans $\sigma$. 
\newline

\para{ Composition}
On définit de manière standard la composition en laissant les deux
stratégies interagir puis en cachant l'interaction dans $B$ (parallel
composition and hiding).
Étant donné deux stratégies $\sigma: A^\Neg \tensor B,\tau :
B^\Neg \tensor C$, on définit la composée
$$
\sigma;\tau = \{ u_{|A,C} \ | \ u \in int(A,B,C) \wedge u_{|A,B} \in
\sigma  \wedge u_{|B,C} \in \tau\} 
$$

De manière classique, on peut montrer la bonne définition de cette
composition en utilisant le lemme fondamental suivant: 

\begin{lemma}[Témoin unique]
  Si $\sigma$ et $\tau$ sont des stratégies de $A^\Neg \tensor B$ et
  $B^\Neg \tensor C$ respectivement, alors pour tout $s\in \sigma;\tau$,
  il existe un unique $u \in int(A,B,C)$ tel que $s =
  u_{|A,C}$ , $u_{|A,B} \in \sigma$ et $u_{|B,C} \in \tau$. 

  De plus, si $s \in \sigma;\tau$ est un préfixe de $ t \in
  \sigma;\tau$, alors le témoin unique de $s$ est préfixe du témoin
  unique de $t$. Ainsi, le lemme du témoin unique s'étend (de manière
  non unique) aux chemins joués par une stratégie.
\end{lemma}

\para{La catégorie des jeux de Conway}
Nous avons maintenant tout ce qu'il faut pour définir une catégorie
compacte close. 

\begin{prop}
La catégorie $\C$, avec pour objets les jeux de Conway et pour
morphismes les stratégies de $A^\Neg \tensor B$, est compacte close

En particulier, on a
$$
\frac{(A \tensor B) \rightarrow C }
     {B \rightarrow A^\Neg \tensor C}
     $$
     
\end{prop}

\para{Trace}
Conséquemment, la catégorie des jeux de Conway est automatiquement
équipée d'une notion de trace. Tout morphisme $f:(X \tensor A)^\Neg
\tensor X \tensor B$ peut être transformé en $\widehat{f}:(X^\Neg
\tensor X)^\Neg \tensor A^\Neg \tensor B$ par commutativité et
associativité du tenseur, et on obtient 

$$Tr_X(f) = \widehat{f}(id_X)$$

\para{$\bool$ool ou Le point du géomètre}

\begin{figure}[h]

\setlength{\columnseprule}{1pt}

\begin{multicols}{2}
  
  \centerline{Partie valide}
  $$
  \begin{array}{ccccc}
    (\bool & \pop & \bool) & \pop & \bool \\
    \vspace{0.2cm}
    
    & & & & \\
    & & & & q \\
    & & q & & \\
      q & & & & \\
      & & & & V \\
      & & & & \\
    \end{array}
  $$
  \centerline{Arène}
  $$
  \xymatrix @-1.5pc {
    && \bool &&\\
    \\
      && \star && \\
    \\
    && \ar@{-}[uu]^{q} . &&\\
    \\
    \ar@{-}[rruu]^{V} . && && \ar@{-}[lluu]_{F} .
  }
  $$
\end{multicols}

\caption{Un exemple de partie valide sans condition de gain} \label{fig:bool1}
\end{figure}
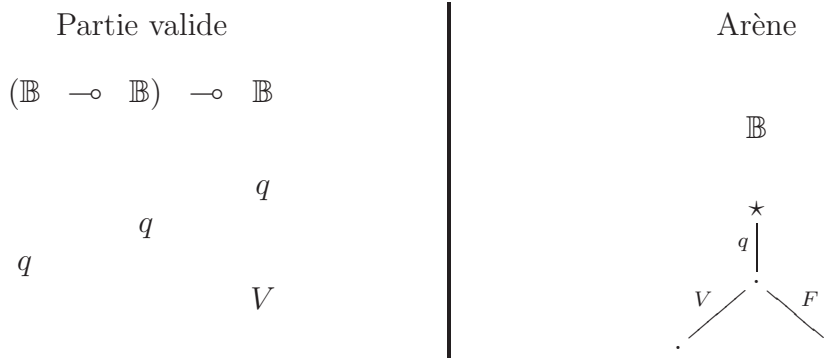

L'exemple \ref{fig:bool1} montre une interaction typique dans le jeu
$(\bool \pop \bool) \pop \bool$ que l'on souhaiterait rejeter. En
effet, ici, le Joueur court-circuite la séquence de questions pour
finalement répondre à la première. Une telle action viole la condition
de bon parenthésage, mais nous n'avons pas pour le moment le
vocabulaire pour interdire un tel comportement.

Il faut donc augmenter le modèle avec une notion de gain qui incorpore
le parenthésage tout en conservant cette belle structure compacte close.

\section*{Jeux de Conway à gain}

Un jeu de Conway à gain $A$ est un n-uplet
$(V_A,E_A,\lambda_A,\kappa_A)$ constitué d'un jeu de Conway
$(V_A,E_A,\lambda_A)$ et d'une fonction de gain dont la structure est
définie ci-dessous.
\newline

\para{La structure du gain}
On définit $\kappa_A = (\kappa_A^+,\kappa_A^-): \paths_A \rightarrow
\mathbb{N}^2$ comme un couple de fonction de $\paths_A$ dans 
$\mathbb{N}$, une pour le gain de Joueur ($\kappa_A^+$), une pour le gain
d'Opposant ($\kappa_A^-$). Il est notable que le gain soit défini pour
tout chemin. On ne peut pas se restreindre aux chemins alternants
car on veut définir le gain d'un jeu tensorisé comme la somme des
gains des projections, qui ne sont pas forcément alternées. On ne peut
pas non plus se restreindre aux parties car on a besoin d'une
condition de gain local pour garantir que toute interaction entre
stratégies gagnantes forment un partie ``bien parenthésée''.   

Nous demandons au gain de vérifier quatre propriétés similaires à
celle demandée à une norme pour un espace vectoriel. 

\begin{description} \label{payoffCond}
\item[compatibilité]
  $
  \quad \forall m \in E_A \quad
  \left \{
  \begin{array}{l}
    \lambda_A(m) = -1 \Implies \kappa_A^+(m) = 0 \\
    \lambda_A(m) = +1 \Implies \kappa_A^-(m) = 0 \\ 
  \end{array}
  \right.
  $
  \\
  Ceci garantit que la distinction Joueur/Opposant a
  bien un sens.  
\item[suffixe dominé] 
  $\quad s:x\chemin y$, $t:y \chemin z$
  $$
  \kappa_A(t) \leq \kappa_A(s;t)
  $$
  Cet axiome exprime qu'une question ne peut pas être répondue dans
  le passé.
\item[sous-additivité] 
  $\quad s:x\chemin y$, $t:y \chemin z$
  $$
  \kappa_A(s;t) \leq
  \kappa_A(s) + \kappa_A(t)   
  $$
  Cet axiome, qui est une sorte d'inégalité de Cauchy-Schwarz,
  stipule qu'un coup ne peut pas poser plus ou moins de questions
  suivant son passé. Ainsi, la composition de deux chemins ne peut que
  faire diminuer le nombre de questions.
\item[norme] 
  $\quad \epsilon_x:x\chemin x \qquad
  \kappa_A(\epsilon_x) = (0,0)
  $
  \\
  C'est la propriété usuelle d'une norme. Dans notre cadre, elle
  exprime qu'en l'absence d'interaction, aucune question ne peut avoir
  été ouverte.
\end{description}

Tous ces axiomes sont assez naturels et ils nous permettent maintenant
de définir une notion de stratégie gagnante dont on va montrer (avec
des arguments algébriques) qu'elle est stable par composition.
\newline

\para{Stratégie gagnante}
Une stratégie est \emph{gagnante} lorsque tout chemin $s$ qu'elle joue
satisfait la condition suivante sur le gain : 

\begin{equation} \label{eqn:payoffCond}
  \kappa_A^+(s) = 0 \Implies \kappa_A^-(s) = 0
\end{equation}

Intuitivement, cette condition exprime qu'une stratégie qui a été
interrogée localement doit ou bien poser une autre question, ou bien
répondre à cette question. Ceci généralise la condition habituelle de
``bon parenthésage''.
\newline

\para{Extension de la structure}
Étant donnés deux jeux $A$ et $B$, on étend les connecteurs logiques
définis sur les jeux de Conway de la manière suivante. 

\begin{description}
\item[Dual]
  Le dual d'un jeu de Conway à gain $A$ est le jeu $A^\Neg = (V_A , E_A ,
  -\lambda_A , (\kappa^-_A,\kappa^+_A))$.
  Le nouveau gain satisfait trivialement les conditions imposées plus
  haut. 
\item[Multiplicatif]
Le produit tensoriel de deux jeux $A$ et $B$ est le produit tensoriel
des jeux de Conway sous-jacents, avec comme gain
$$
\forall s \in \paths_{A \tensor B} \qquad  \kappa_{A \tensor B}(s) =
\kappa_A(s_{|A}) + \kappa_B(s_{|B}) 
$$ 
\end{description}

Pour s'assurer que ce nouveau gain satisfait aux conditions demandées,
nous avons besoin de remarquer les propriétés suivantes:

\begin{enumerate}
\item
  Comme tout coup de $A \tensor B$ appartient soit à $A$, soit à  $B$,
    la propriété de compatibilité s'étend naturellement.
\item
  Étant donnés deux chemins $s:x\path y$ et $t:y \path z$, on a
  \begin{eqnarray*}
    \kappa_{A\tensor B}(t) & = & \kappa_A(t_{|A}) + \kappa_B(t_{|B})\\
    & \leq & \kappa_A((s;t)_{|A}) + \kappa_B((s;t)_{|B})\\
    & \leq & \kappa_{A\tensor B}(s;t)
    \end{eqnarray*}
\item   
  Étant donnés deux chemins $s:x\path y$ et $t:y \path z$, on a
  \begin{eqnarray*}
      \kappa_{A\tensor B}(s;t) & = & \kappa_A((s;t)_{|A}) + \kappa_B((s;t)_{|B})\\
      & \leq & (\kappa_A(s_{|A}) + \kappa_A(t_{|A})) +
      (\kappa_B(s_{|B}) + \kappa_B(t_{|B})) \\
      & \leq & \kappa_{A\tensor B}(s) + \kappa_{A\tensor B}(t)
  \end{eqnarray*}
\item
  $\forall x,y \quad \kappa_A(\epsilon_{x\tensor y }) = \kappa_A(\epsilon_x) +
  \kappa_A(\epsilon_y) = (0,0)
  $
\end{enumerate}

Il est urgent de vérifier que la condition de gain sur stratégie est
préservée par composition. C'est toujours un point délicat en
sémantique des jeux, et nous espérons que la formulation du gain sous
forme d'axiomes permet de rendre la démonstration suivante un peu
moins indigeste pour le lecteur.

\begin{prop}
  Soit deux stratégies gagnantes $\sigma:A ^\Neg \tensor B$ and $\tau:B
  ^\Neg \tensor C$. 
  La stratégie $\sigma;\tau:A^\Neg \tensor C$ est aussi gagnante. 
\end{prop}

\begin{proof}
  
  Nous procédons par l'absurde en supposant qu'il existe un chemin $s$
  (présupposé le plus petit) joué par la stratégie $\sigma;\tau$
  tel que $\kappa_{A ^\Neg \tensor C}(s) \in 0 \times \mathbb{N}^*$.
  On montre alors que soit $\sigma$, soit $\tau$ a triché.
  
  Dans un premier temps, en utilisant le lemme du témoin unique sur
  $s$, on obtient un $u$ tel que $u _{|A,B}$ est joué par $\sigma$, $u
  _{|B,C}$ est joué par $\tau$ et $u$ est dans $int(A,B,C)$.
  
  Alors, en utilisant la \emph{domination par suffixe}, on a pour
  tout $u'$ préfixe de $u$, $\kappa^+_{A^\Neg \tensor C}(u')=~0$, ce
  qui entraîne en utilisant la définitions du produit tensoriel  
  $$\kappa_{A^\Neg}^+(u'_{|A}) = \kappa_C^+(u'_{|C}) = 0$$   
  
  À ce moment, deux cas doivent être considérés
  \begin{enumerate}
    \item
      ${\kappa_B(u_{|B}) = (0,0)}$. \\
      Alors $\kappa_{A^\Neg \tensor B}(u) = \kappa_{A^\Neg}(u)$ et
      $\kappa_{B^\Neg \tensor C}(u) = \kappa_C(u)$. 
      Mais comme $\kappa_{A^\Neg \tensor C}(u)^- =  \kappa_{A^\Neg}^-(u) +
      \kappa_C(u)^- > 0 $ par hypothèse, on a que ou bien $
      \kappa_{A^\Neg}^-(u) > 0$, ou bien $\kappa_{C}^-(u) >0$.
      
      Ainsi, soit $\kappa_{A^\Neg \tensor B}(u)\in 0 \times \mathbb{N}^* $, soit
      $\kappa_{B^\Neg \tensor C}(u)\in 0 \times \mathbb{N}^* $, ce qui
      implique qu'au moins une des deux stratégies a triché.
      
    \item
      ${\kappa_B(u_{|B}) \neq (0,0)}$. \\
      Considérons maintenant $v$, le plus petit suffixe de $u$ tel que
      ${\kappa_B(v_{|B}) \neq (0,0)}$.

      Soit $v = m;v' $.
      Comme $v$ est le plus suffixe à gain non nul, nécessairement
      $\kappa_B(v'_{|B}) = (0,0)$ et alors $\kappa_B(v_{|B}) \leq
      \kappa_{B}(m)$ (par \emph{sous-additivité}).
      Mais comme $u$ est la plus petite interaction donnant lieu à un
      mauvais gain dans $A^\Neg \tensor C$, $m$ est nécessairement
      dans $E_B$. Par la condition de \emph{compatibilité}, on sait que
      soit $\kappa_B^+(m) = 0$, soit $\kappa_{B}^-(m) = 0$, ce qui
      implique que ou bien $\kappa_B^+(v_{|B}) = 0$, ou bien
      $\kappa_B^-(v_{|B}) = \kappa_{B^\Neg}^+(v_{|B}) = 0$.

      Traitons le cas $\kappa_{B^\Neg}^+(v_{|B}) = 0$ (l'autre étant
      identique si on remplace $\tau$ par $\sigma$) ce qui entraîne
      que $\kappa_{B^\Neg}^-(v_{|B}) > 0$. 
      Dans ce cas, considérons le chemin parcouru par $\tau$.
      Des considérations ultérieures montrent que $\kappa_{B^\Neg
      \tensor C}^+(v_{|B^\Neg \tensor C}) = 0$. Mais
      $\kappa_{B^\Neg}^-(v_{|B}) > 0 \Implies \kappa_{B^\Neg \tensor
      C}^-(v_{|B^\Neg \tensor C}) > 0$ par définition du gain sur le
      produit tensoriel. Alors la stratégie $\tau$ à joué un chemin interdit.
  \end{enumerate}
  
\end{proof}

\para{La catégorie $\p$ des jeux de Conway à gain}
La catégorie $\p$, qui a pour objets les jeux de Conway à gain et pour
morphismes les stratégies gagnantes entre ces jeux, est une catégorie
compacte close. 
En effet, la notion de gain surajoutée n'a pas perturbé la structure
auto-duale du tenseur (car la somme est commutative dans $\mathbb{N}$).
\newline

\para{$\bool$ool ou Le point du géomètre (2)}

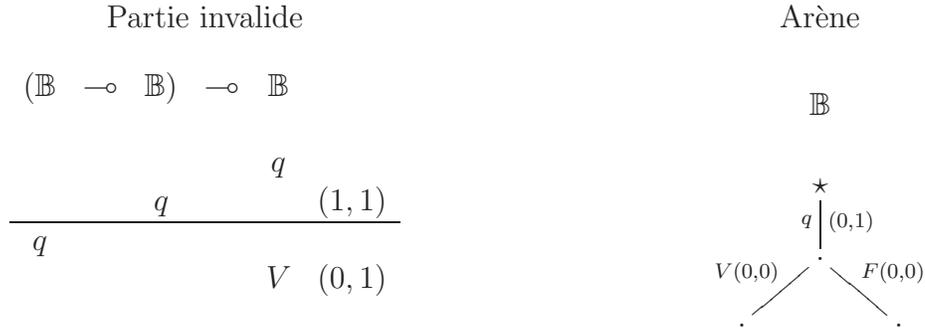
\begin{figure}[h]

\begin{multicols}{2}

  \setlength{\columnseprule}{1pt}

  \centerline{Partie invalide}
  $$
  \begin{array}{cccccc}
    (\bool & \pop & \bool) & \pop & \bool & \\
    & & & & &\\
    & & & & q & \\
    & & q & & & (1,1)\\
    \hline
    q & & & & & \\
    & & & & V & (0,1)\\
    & & & & & \\
  \end{array}  
  $$
  \centerline{Arène}
  $$
  \xymatrix @-1.5pc {
    && \bool &&\\
    \\
    && \star && \\
    \\
    && \ar@{-}[uu]^{q}_{(0,1)} . &&\\
    \\
    \ar@{-}[rruu]^{V (0,0)} . && && \ar@{-}[lluu]_{F (0,0)} .
  }
  $$
\end{multicols}

\caption{L'exemple \ref{fig:bool1} est maintenant interdit par la
  condition de gain} \label{fig:bool2}

\end{figure}

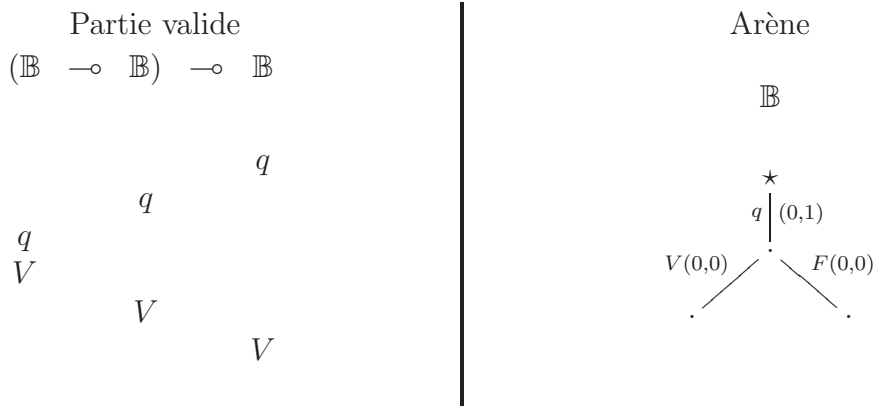
\begin{figure}[h]

\setlength{\columnseprule}{1pt}

\begin{multicols}{2}
    \center Partie valide
    $$
    \begin{array}{cccccc}
      (\bool & \pop & \bool) & \pop & \bool & \\
      \vspace{0.2cm}
      
      & & & & &\\
      & & & & q & \\
      & & q & & & \\
      q & & & & & \\
      V & & & & & \\
      & & V & & & \\
      & & & & V & \\
      & & & & & \\
    \end{array}
    $$
    \center Arène
    $$
    \xymatrix @-1.5pc {
      && \bool &&\\
      \\
      && \star && \\
      \\
      && \ar@{-}[uu]^{q}_{(0,1)} . &&\\
      \\
      \ar@{-}[rruu]^{V (0,0)} . && && \ar@{-}[lluu]_{F (0,0)} .
    }
    $$
\end{multicols}

\caption{Exemple d'une interaction totale} \label{fig:bool3}

\end{figure}

Revenons sur l'exemple \ref{fig:bool2}. On peut maintenant exprimer
que le premier coup Opposant du jeu booléen est une question et que
Joueur y répond par Vrai ou Faux. Ainsi, l'interaction précédente est
maintenant non valide (voire figure \ref{fig:bool2}).
Le Joueur est maintenant obligé de respecter la séquence de questions
ouvertes dans l'interaction (voire figure \ref{fig:bool3}. 

C'est la condition de bon parenthésage. 
\newline

\para{Une revisite du parenthésage}
Nous décrivons maintenant la notion de parenthésage induite par le
gain :

\begin{itemize}
\item
  Une partie $s$ est bien parenthésée pour Joueur
  lorsque tout chemin de longueur paire $t$ dans $s$ qui se termine
  par un coup Joueur vérifie $\kappa_A^+(t) = 0 \Implies
  \kappa_A^-(t) = 0$ 
\item
  Une partie $s$ est bien parenthésée pour Opposant
  lorsque tout chemin de longueur paire $t$ dans $s$ qui se termine
  par un coup Opposant vérifie $\kappa_A^-(t) = 0 \Implies
  \kappa_A^+(t) = 0$ 
\item
  On dit qu'une partie $s$ est \emph{bien parenthésée} lorsque qu'elle
  est bien parenthésée à la fois pour Joueur et pour Opposant
\end{itemize}

Nous voulons montrer que l'interaction entre deux stratégies gagnantes
fournit toujours une partie bien parenthésée.
Pour cela, il nous faut d'abord donner un sens aux mots interaction
entre stratégies.
Nous avons besoin du jeu $\point$ qui correspond au graphe $\bullet
\xrightarrow{o} \bullet$ qui contient un unique coup Opposant de gain
nul (on note $\point$ par analogie avec la catégorie $\point$ définie
par Lawvere). 

\begin{definition}[interaction entre stratégies]
Soit $\sigma : A$ et $\tau : A^\Neg \tensor \point$.

Intuitivement, on définit $\sigma \Join \tau$ comme la partie où le
premier coup est la réponse de $\tau$ au coup Opposant dans $\point$
et où le reste est déduit des réponses respectives de $\sigma$ et
$\tau$. 

Formellement, l'interaction est définie par
$$
\sigma \Join\tau = \{\epsilon\} \cup \{s \cdot m \ | \ s \cdot m \in
\sigma \wedge o \cdot s \in \tau \}
$$
où $s$ est une partie, $m$ un coup joueur et $o$ l'unique coup de
$\point$.
\end{definition}

On déduit de ces définitions deux propriétés évidentes sur les
stratégies gagnantes. 

\begin{prop}
{\bf (Bon parenthésage)} \\
\begin{enumerate}
\item
  Toute stratégie gagnante joue des parties bien parenthésées pour Joueur.
\item
  Toute interaction entre deux stratégies gagnantes $\sigma:A$ and
  $\tau:A^\Neg \tensor \point$ produit une partie bien parenthésée.
\end{enumerate}
\end{prop}

\begin{proof}
  \begin{enumerate}
    \item
      La première propriété est une conséquence directe de la
      définition de stratégie gagnante. Ceci n'est pas surprenant car
      c'est la notion de parenthésage qui à motiver notre
      définition du gain est stratégie gagnante dans les jeux de
      Conway.
    \item
      Cette deuxième propriété mérite plus d'attention.
      Considérons un chemin de longueur paire $s \cdot m$ dans
      l'interaction $\sigma \Join \tau$.
      \begin{description}
	\item[Premier cas: $m$ est un coup joueur.]
	  Alors $s \cdot m \in \sigma$, ce qui implique par la
	  propriété ci dessus que $\kappa_A^+(s \cdot m) = 0 \Implies
	  \kappa_A^-(s \cdot m) = 0$.  
	\item[Deuxième cas: $m$ est un coup opposant.]
	  Alors $o \cdot s \cdot m \in \tau$, ce qui implique par la
	  propriété ci dessus que $\kappa_{A^\Neg \tensor \point}^+(o
	  \cdot s \cdot m) = 0 \Implies \kappa_{A^\Neg \tensor
	  \point}^-(o \cdot s \cdot m) = 0$.
	  Or, en utilisant la définition du gain sur le tenseur et
	  le dual (et aussi que le gain dans $\point$ est nul), on
	  obtient la propriété: $\kappa_{A}^-(s \cdot m) = 0 \Implies
	  \kappa_{A}^+(s \cdot m) = 0$  
      \end{description}
  \end{enumerate}
\end{proof}

\para{Bilan}
Nous avons à présent à notre disposition une catégorie compacte close
avec une notion de contrôle induite par un gain algébrique. Cette
catégorie possède de plus le bon goût de rester dans notre cadre
asynchrone et nous pouvons ainsi espérer, dans un futur proche, y
introduire la notion d'homotopie.  
\newline

Fort de ce premier pas encourageant, il nous faut maintenant nous
tourner vers la construction algébrique de l'exponentielle.

\chapter[Construction du comonoïde commutatif libre]{La construction du comonoïde commutatif libre dans les jeux
  de Conway} \label{ch:exponentielle}

La modalité exponentielle est à priori le point faible de la logique
linéaire du point de vue structurelle. En effet, c'est la seule
construction qui, étant donné un modèle, n'apparaît pas comme
unique. En effet, contrairement aux autres connecteurs, pas moyen de
prouver l'équivalence (en terme de théorie de la preuve) de deux
modalités satisfaisant aux règles de l'exponentielle.
Pourtant du point de vue logique, elle est au c\oe ur de la
magie linéaire et de sa capacité à décortiquer à peu près n'importe
quel formalisme logique. Il paraît donc important de comprendre
comment réconcilier structure catégorique et règles logiques en ce
point.
Dans ce qui suit, nous motivons la recherche du comonoïde commutatif
libre même si nous présentons par la suite la construction du monoïde
commutatif libre.
Cette construction duale plus canonique en catégorie a eu nos faveurs
pour une présentation générale. Le lecteur n'aura qu'à inverser
quelques flèches pour retrouver la construction de l'exponentielle.

\section*{Le chemin de la liberté}

Dans un premier temps, regardons ce qui se passe dans un modèle
où l'on peut définir plusieurs modalités exponentielles, les fameux
espaces de cohérence.
\newline

\para{Le modèle relationnel}
Plus besoin de rappeler cette catégorie \emph{tarte à la crème}
(le lecteur peut se reporter à \cite{Girard:linear-logic} s'il n'est
pas familier avec ce mets), tant ce modèle est à la base de nombreuses
intuitions sémantiques. Ici encore, il va nous permettre de comprendre
les différences que peuvent présenter deux modalités exponentielles au
sein d'un même modèle.
Historiquement, la première modalité due à Girard est définie comme
l'ensemble des cliques finies, alors que la modalité la plus usitée
maintenant est celle des multicliques finies. Pourquoi cette
préférence envers les multi-ensembles? Il s'avère que ces derniers
donnent lieu à une construction libre alors que la modalité avec
cliques finies ne donne lieu qu'à un comonoïde vérifiant les
propriétés demandées par la logique linéaire.

Il semble donc que la construction libre soit une bonne approche à ce
problème. En effet, comme tout objet libre, le comonoïde devient alors
unique modulo isomorphisme. Il retrouve alors la qualité des
constructions structurelles comme le produit tensoriel.
\newline

\para{Les catégories de Lafont}
Dans sa revisite des modèles catégoriques de la logique linéaire
\cite{mellies:categorical-models}, Paul-André Melliès rappelle la
définition de Lafont d'une catégorie symétrique monoïdale fermée avec
produits finis et comonoïde commutatif libre sur chaque objet de la
catégorie. Même si ce formalisme ne capture pas certains modèles comme
celui susmentionné des espaces de cohérence avec cliques finies, il
présente l'avantage d'être un cas pas si particulier et très simple
des catégories linéaires (résultat dû à Bierman). Nous voilà donc avec
un deuxième argument pour rechercher des comonoïdes commutatifs libres
dans une catégorie sensée capturer la logique linéaire. Reste
maintenant à définir un cadre pour en assurer l'existence et le
calculer. 

\section*{La boîte à outils du catégoricien}

\para{Attention}
Dans un souci de clarté et de simplicité, nous considérons dans tout
ce qui suit que toutes les colimites mentionnées existent. Cela évite
de supposer la catégorie cocomplète, ce qui est très fort (et faux
dans les jeux de Conway), ou de mentionner à chaque fois que l'on
suppose que la colimite existe, ce qui est très lourd.
\newline

Comme annoncé ci-dessus, nous allons maintenant décrire des pistes de
construction du monoïde commutatif libre s'appuyant sur des travaux en théorie des
catégories. Dans un premier temps, nous présentons nos travaux pour
recomprendre cette construction dans le cas connu où le tenseur
commute aux colimites. Nous présentons ensuite un résultat méconnu de
Dubuc~\cite{dubuc74} qui permet d'étendre cette construction pour le
monoïde libre dans le cas où le
tenseur commute aux colimites $\varphi$-filtrées (voir annexe).
Nous avons donc dû adapter ce résultat pour retrouver le monoïde
commutatif libre.
Enfin, nous aborderons une possible extension vers une théorie
monoïdale à l'instar des théories algébriques de Lawvere.
\newline

\para{La revisite d'un cas classique avec les extensions de Kan}
Dans le cas où le tenseur commute aux colimites, il est bien connu que
le monoïde commutatif libre se construit par la formule 
$$\Sigma A = \bigoplus_n S^n(A)$$
où $S^n(A)$ est le symétrisé de $A^{\otimes n}$.
Pour recomprendre cette formule, considérons l'extension de Kan à
gauche (voir annexe) sur
$\mathcal{S}_0$ (la catégorie des ensembles finis et fonctions
ensemblistes) du jeu $A$ vu comme un foncteur monoïdal de
$\mathbb{B}$ (la catégorie des ensembles finis et 
fonctions bijectives) dans la catégorie $\C$ des jeux de Conway. 

$$
\xymatrix{
\C & & \\
& & \\
\mathbb{B} \ar[uu]^{A} \ar[rr]^{J} & & \mathcal{S}_0 \ar[lluu]_{\exists_J(A)}
}
$$

Les travaux de Day et Street \cite{day95} sur la monoïdalité de
l'extension de Kan avec la convolution comme structure monoïdale sur
la catégorie des foncteurs (voir annexe) nous donne l'équation
suivante:

$$
\exists_J(A *_\mathbb{B} B) \cong \exists_J(A) *_{\mathcal{S}_0}
\exists_J(B) 
$$

Mais, de manière très agréable, lorsque le tenseur commute avec les
colimites, on obtient

\begin{prop}\label{prop:convolution}
Si la catégorie monoïdale ($\C$,$\tensor$,$I$) à un tenseur qui
commute aux colimites, alors
$$
\left \{
\begin{array}{rcl}
A *_\mathbb{B}B (n) & \cong & S^{n}(A \oplus B) \\
A *_{\mathcal{S}_0}B (n) & \cong & S^{n}(A) \tensor S^{n}(B)
\end{array}
\right.
$$
\end{prop}

\begin{proof}
  \begin{eqnarray*}
  A *_\mathbb{B} B (n) & \equiv & \int^{m,m'} \mathbb{B}(m+m',n).(A^m
  \tensor B^{m'}) \\
  & \cong & \bigoplus_{m} S(A^m \tensor B^{n-m}) \\
  & \cong & S^n(A \oplus B) \\
  & & \\
  A *_{\mathcal{S}_0} B (n) & \equiv & \int^{m,m'} \mathcal{S}_0(m+m',n).(A^m
  \tensor B^{m'}) \\
  & \cong &  \int^{m,m'} (\mathcal{S}_0(m,n) \times \mathcal{S}_0(m',n)).(A^m
  \tensor B^{m'}) \\
  & \cong &  (\int^{m} \mathcal{S}_0(m,n).A^m) \tensor (\int^{m}
  \mathcal{S}_0(m,n).B^m) \\
  & \cong & S^n(A) \tensor S^n(B)
  \end{eqnarray*}
\end{proof}

Ce qui donne donc lorsque applique en $1$ (si on note $\exists_J(A)(1)
\equiv {\Sigma A}$)
\begin{equation} \label{eq:extensionMonoidale}
{\Sigma (A \oplus B)} \cong {\Sigma A} \tensor {\Sigma B}
\end{equation}

On en déduit alors que $\Sigma A$ est le monoïde commutatif libre sur
$A$. De plus, la construction canonique des extensions de Kan à gauche
nous redonne la formule 
$$ 
\Sigma A \equiv \bigoplus_n S^n(A)
$$
Regardons maintenant comme illustration, une alternative aux jeux de
Conway, les jeux de Conway synchronisés.
\newline

\para{Jeux de Conway synchronisés}
Malheureusement, le tenseur des jeux de Conway ne commute pas avec la
somme. Pour arranger cela (momentanément car ce n'est pas le modèle
qui nous intéresse pour le moment), il faut définir une version
synchronisée du tenseur. Nous restons ici dans le domaine de
l'intuition car la formalisation ne nous parait pas nécessaire à la
compréhension.

Lorsque l'on regarde si le tenseur commute avec la somme dans les jeux
de Conway positifs, on s'aperçoit que l'on ne peut pas décrire la
stratégie 
$$
A \tensor (B \oplus C) \longrightarrow (A \tensor B) \oplus (A \tensor C)
$$ 
En effet, lorsque Opposant joue (à gauche) son premier coup dans $A$,
on ne sait pas quel jeu entre $B$ et $C$ ``sacrifier'' (à droite). Pour
résoudre ce problème, il faut passer au tenseur synchronisé
$\tensor_s$ qui demande à ce que le premier coup soit joué dans les
deux composantes en même temps. Ainsi, on obtient l'équation
$$
A \tensor_s (B \oplus C) \cong (A \tensor_s B) \oplus (A \tensor_s C)
$$ 
On en déduit donc dans ce cadre que le monoïde libre se calcule
simplement par la formule susmentionnée. Cela exprime que dans le
cadre synchronisé, Opposant doit annoncer le nombre de copies qu'il va
jouer le reste de la partie. On a donc bien la stratégie 
$$
\begin{array}{ccccc}
  (\Sigma A  & \tensor_s & \Sigma A) & \longrightarrow & \Sigma A \\
  & & & & \\
  n & & m & &  \\
  & & & & n+ m \\
\end{array}
$$
On ne peut en revanche pas faire la même chose si le tenseur n'est pas
synchronisé car dans ce cas Joueur ne sait pas combien de copies ouvrir
(à droite). 

C'est donc bien la même raison qui brise la structure de
monoïde pour $\Sigma A$ et qui empêche le tenseur de commuter avec la somme. 
Voyons comment s'affranchir de ce problème.
\newline

\para{Une revisite de la construction de Dubuc}
Il s'avère que le monoïde libre peut s'obtenir de façon similaire
lorsque le tenseur commute uniquement aux colimites
$\varphi$-filtrées~\cite{dubuc74}. Pour cela, il faut bien évidemment
étendre la catégorie sur laquelle on travaille en une catégorie
$\varphi$-filtrée. Dubuc le réalise en définissant une sorte de catégorie
simpliciale étendue aux ordinaux (plus un point particulier $-1$)
$\Delta_{ord}$ qui possède l'agréable propriété d'être $\varphi$-filtrée
pour tout ordinal limite $\varphi$.

{\bf Remarque.}
Les constructions suivantes s'appliquent uniquement aux objets pointés
(ie. avec une flèche $\mu_A:I \rightarrow A$)
et au morphisme entre tels objets (ie. $f:A \rightarrow B$ telle que
$\mu_A;f =\mu_B$). On a donc pas le monoïde commutatif libre général,
mais une version restreinte aux objets pointés.
Cela ne pose pas de problème pour les jeux de
Conway car tout monoïde à une unique flèche de
$I$ dans $A$.
\newline 

Cette construction est malheureusement assez technique lorsque l'on
regarde les détails et nous donnons ici une revisite simple de cette
construction lorsque $\varphi = \omega$. 

Dans ce cas, il suffit de considérer la catégorie simpliciale usuelle
$\Delta$ des entiers et fonctions croissantes. On procède ensuite
comme suit

\begin{itemize}
\item
  On prend un objet $A$ dans $\C$ pointé au sens où il y a une flèche
  $I \rightarrow A$
\item 
  On construit un foncteur monoïdale $T$ de $\Delta$ dans $\C$
  envoyant $0$ en $I$, $1$ en $A$ et l'unique flèche $0 \rightarrow 1$
  dans $I \rightarrow A$, vérifiant donc
  $$
  T(m + n) = T(m) \tensor T(n)
  $$
  pour $m,n$ entiers
\item
  on pose 
  $$
  TA = \underrightarrow{colim}(\Delta \xrightarrow{T} \C)
  $$
\end{itemize}

Définissons maintenant
$$
\widetilde T : 
\left \{
\begin{array}{rcl}
  \Delta\times\Delta & \rightarrow& \C \\
  (m,n) & \mapsto & T(m + n) = T(m) \tensor T(n) 
\end{array}
\right.
$$
En utilisant que le tenseur commute aux colimites filtrées, on obtient:
$$
TA \tensor TA = \underrightarrow{colim}(\widetilde T)
$$ 
Or, il est clair que le cône de $TA$ sur $T$ s'étend en un cône sur
$\widetilde T$. 
L'universalité de $TA \tensor TA$ dans la catégorie des cônes sur
$\widetilde T$ nous fournit la flèche
$$
TA \tensor TA \rightarrow TA
$$
qui, avec la flèche $1 \rightarrow TA$ venant du cône, font de $TA$ un
monoïde dans $\C$ (il faut encore vérifier les deux diagrammes d'un
monoïde, diagrammes qui passent tous seuls en utilisant qu'une
colimite est un objet initial dans la catégorie des cônes).

On peut maintenant statuer sur le théorème suivant.

\begin{theoreme}
  Le monoïde $TA$ est le monoïde libre sur l'objet pointé $A$   
\end{theoreme}

Afin d'obtenir la liberté du monoïde, il nous faut prouver
le lemme suivant

\begin{lemma}
  Soit $M$ un monoïde avec une flèche entre objets pointés $f:A \rightarrow
  M$. 
  Alors, on peut construire un cône $\cone{C}_M$ sur $T$.
\end{lemma}  

\begin{proof}
  Les flèches de $1$ et $A$ dans $M$ sont déjà données (rappelons que
  $M$ est un monoïde). La commutation du diagramme
  $$
  \xymatrix{
    1 \ar[d] \ar[r] & A \ar[ld] \\
    M &
  }
  $$
  est précisément la propriété requise pour être une flèche entre objets
  pointés.
  Reste à construire les flèches de $A^n$ dans $M$. Rien de plus
  simple avec la multiplication $d$
  $$
  \xymatrix{
    A^n \ar[r]^{f^n} & M^n \ar[r]^{d^n} & M
  }
  $$
  Comme le monoïde n'est pas commutatif, l'ordre d'application de $d$
  importe, et nous choisissons une ordre par la gauche.
  Il est ensuite évident grâce aux propriétés d'un monoïde que le
  diagramme
  $$
  \xymatrix{
    A^m \ar[d]_{f^m} \ar[r]^{\pi} & A^n \ar[d]^{f^n} \\
    M^m \ar[d]_{d^m} & M^n \ar[ld]^{d^n} \ar[l]^{d^{n-m}}\\
    M &
  }
  $$
  commute pour tout $m\leq n$
\end{proof}

Considérons maintenant un monoïde $M$ avec une flèche entre objets
pointés $f:A \rightarrow M$. Il existe donc un cône $\cone{C}_M$ sur
$T$. On en déduit qu'il existe une flèche entre $TA$ et $M$ qui fait
commuter
$$
\xymatrix{
  A \ar[d] \ar[r] & TA \ar[ld] \\
  M &
}
$$
dont il
nous faut montrer qu'elle est monoïdale.
    
De la même manière que dans le lemme précédent, on construit un cône
de $M$ sur $\widetilde T$ et un cône de $M \tensor M$ sur $\widetilde
T$. On en déduit alors que le diagramme
$$
\xymatrix{
  1 \ar@{=}[d] \ar[r] & TA \ar[d] & TA \tensor TA \ar[l] \ar[d] \\
  1 \ar[r] & M & M \tensor M \ar[l]   
}
$$
commute, le seul point délicat étant le carré de droite dont on montre
que les deux flèches sont égales car elles factorisent toutes les deux
le cône de $M$ sur $\widetilde T$.  

La flèche entre $TA$ et $M$ venant de la propriété de colimite est
donc monoïdale.
Reste à montrer qu'elle est unique. Mais étant donné une flèche
monoïdale entre $TA$ et $M$, il est facile de montrer que celle-ci
factorise $\cone{C}_M$, elle est donc unique. 
\newline

\para{Extension au monoïde commutatif libre}
Nous devons maintenant adapter le résultat de Dubuc pour pouvoir
décrire l'exponentielle. À cette fin, nous changeons de catégorie de
base en passant de la catégorie simpliciale à la catégorie des
injections. Ainsi, on fait apparaître les permutations de $n$ dans $n$
représentant la commutation. On note $Inj$ cette catégorie. On suppose
toujours que notre catégorie $\C$ commute aux
colimites $\omega$-filtrées. Dans ce cadre, on peut 
étendre sans difficultés le foncteur $T : \Delta \rightarrow \C$
à un foncteur $T' : Inj \rightarrow \C$ pour un objet pointé $A$.

Définissons par anticipation
$$
\Sigma A = \underrightarrow{colim}(\Delta_\omega
\xrightarrow{T'} \C)
$$

On va maintenant montrer la propriété suivante

\begin{prop}\label{prop:monCom}
  $\Sigma A$ est le monoïde commutatif libre sur $A$ 
\end{prop}

Pour la démonstration, nous avons besoin de pouvoir étendre le cône
$\cone{C}_M$ sur $T$ en un cône $\cone{C}'_M$ sur $T'$ lorsque l'objet
sous-jacent est commutatif.

\begin{lemma}
  Soit $M$ un monoïde commutatif une flèche entre objets pointés $f:A
  \rightarrow M$.
  Alors le cône $\cone{C}_M$ sur $T$ peut être étendu en un cône
  $\mathfrak{C}'_M$ sur $T'$.  
\end{lemma}

\begin{proof}
Comme les domaines de $T$ et $T'$ sont identiques, il suffit de
montrer que les permutations sont préservées par le cône (les autres
flèches sont ensuite déduites par composition).
Notre objet $M$ étant un monoïde commutatif, on peut donc montrer
qu'étant donner une permutation $\pi : A^n \rightarrow A^n$, le
diagramme 
$$
\xymatrix{
  A^n \ar[r]^{\pi} \ar[d]_{f^n} & A^n \ar[d]^{f^n} \\
  M^n \ar[r]^{\alpha} \ar[d]^{d}& M^n \ar[dl]^{d}\\
  M & \\
}
$$ 
commute, et ainsi $\mathfrak{C}_M$ s'étend naturellement en un cône
$\mathfrak{C'}_M$ sur $T'$.
\end{proof}

\begin{proof}[Preuve de la proposition \ref{prop:monCom}]
  Soit un monoïde commutatif $M$ avec une flèche d'objet pointé $f:A
  \rightarrow M$. 
  D'après le lemme ci-dessus, $M$ possède un cône sur $T'$. On
  en déduit une flèche $\Sigma f:\Sigma A \rightarrow M$.
  Il faut maintenant montrer qu'elle est monoïdale 
  Comme le tenseur préserve les colimites filtrées, on a les égalités
  $$
  \left \{
  \begin{array}{rcl}
    TA \tensor TA & = &
    \underrightarrow{colim}(\widetilde T) \\
    \Sigma A \tensor \Sigma A & = &
    \underrightarrow{colim}(\widetilde{T'})
  \end{array}
  \right.
  $$
  D'après les lemmes précédents, on obtient les cônes
  $\mathfrak{C}_M$ et $\mathfrak{C'}_M$ de $M$ sur la catégorie
  simpliciale et injective, ainsi que les cônes de $M\tensor M$ sur
  les catégories produits.   
  On peut donc construire le diagramme
  $$
  \xymatrix{
    1 \ar@{=}[d] \ar[rr] & & TA \ar[d] \ar@/^7.5ex/[dd]
    & & TA \tensor TA \ar[ll] \ar[d] \ar@/^7.5ex/[dd]\\ 
    1 \ar@{=}[d] \ar[rr] & &  \Sigma A \ar[d] & & \Sigma A \tensor \ar'[l][ll]
    \ar[d] \Sigma A\\ 
    1 \ar[rr] & & M & & M \tensor M \ar[ll]   
  }
  $$
  dont on déduit la commutativité de l'unicité des flèches.
  On a montré qu'il y avait une flèche monoïdale entre
  $\Sigma A$ et $M$. L'unicité vient encore du fait qu'une flèche
  monoïdale entre ces deux objets factorise le cône sur $M$.

  $\Sigma A$ est bien le monoïde commutatif libre.
\end{proof}

\para{Remarque}
Il semble que cette construction peut se généraliser aux catégories
dont le tenseur ne commute qu'aux colimites $\varphi$-filtrées mais la
présentation de ce résultat nécessite un cadre trop lourd à mettre en
place pour pouvoir apparaître dans ce rapport. 

\section*{Le monoïde commutatif libre vu comme une $\Sigma$-algèbre libre}

Nous mentionnons brièvement ici la structure monadique qui découle de
la fonction sur les objets $\Sigma : A \rightarrow \Sigma A$. En
effet, on peut montrer qu'étant donnée une flèche entre deux objets
pointés $f : A \rightarrow B$, on peut construire une flèche $\Sigma
f: \Sigma A \rightarrow \Sigma B$ qui fait de $\Sigma$ un
endofoncteur dans la catégorie $\C$.

On sait déjà qu'il existe une flèche $\mu: A \rightarrow \Sigma A$,
mais on a aussi le diagramme suivant (grâce à la liberté)
$$
\xymatrix{
  \Sigma A \ar[d]  \ar@/^6ex/[dd]^{id}\\
  \Sigma A \tensor \Sigma A \ar[d] \\
  \Sigma A
}
$$
On en déduit que $\Sigma A \trianglelefteq \Sigma \Sigma A$ (rétraction).

On vérifie alors aisément que $\Sigma$ est une monade sur $\C$, ce qui
permet de voir les monoïdes libres comme des $\Sigma$-algèbres libres. 

\section*{Construction dans les jeux de Conway}

\para{Attention}
Ce qui suit utilise les résultats précédents en prenant à chaque fois
les notions duales. Au lieu de construire le monoïde commutatif
libre, on construit le comonoïde commutatif libre; au lieu de regarder
les colimites filtrées, on regarde les limites filtrées; etc \ldots 
\newline

\para{Exponentielle}
Dans un premier temps, résumons ce dont nous avons besoin pour
construire le comonoïde (commutatif) libre sur un objet $A$ à l'aide
de la construction de Dubuc. Il faut 
\begin{itemize}
\item
  que notre objet $A$ soit pointé au sens où il possède une flèche $A
  \rightarrow I$
\item
  que la catégorie dans laquelle on travaille ait son tenseur qui
  commute aux limites $\omega$-filtrées pour un certain $\omega$.
\item
  que toutes les limites considérées existent bien dans la catégorie.
\end{itemize} 

Comme cela est possible, nous allons construire le comonoïde libre
pour n'importe quel objet de $\C$.

Nous partons d'abord d'un objet $A \in \n$ de la catégorie des jeux de
Conway négatif . Celui-ci est évidemment
pointé car l'ensemble des parties commençant par Opposant de $A
\pop I$ est vide. La stratégie canonique entre $A$ et $I$ est donc la
stratégie vide.
De plus, comme le tenseur à un adjoint à gauche dans $\p$, il commute
aux limites et en particulier aux limites $\omega$-filtrées. On sait
donc que si la limite du diagramme suivant existe, c'est notre
comonoïde commutatif libre 
\vspace{0.3cm}
$$
\xymatrix{
1  & A \ar[l] & A^{\tensor 2} \ar@(ul,ur)[] \ar@2[l] & A^{\tensor 3}
\ar@(ul,ur)[] \ar@3[l] & \ldots
}
$$

Cette limite se calcule et on obtient :

$$
{!A} = S^\infty(A) 
$$
où $S^\infty(A)$ signifie le tenseur infini symétrisé de $A$ que l'on
peut décrire comme le jeu ${!A} = (V_{!A},E_{!A},\lambda_{!A})$ 

\begin{itemize}
\item
  $V_{!A}$ est l'ensemble des fonctions de $\mathbb{N}$ dans $V_A$ qui
  envoient $i$ premiers entiers sur une position différente de
  $\star_A$ et le reste sur $\star_A$ (pour un certain $i$).
\item
  $\star_{!A}$ est la fonction constante qui envoie tout le monde sur
  $\star_A$ 
\item
  Il existe un coup entre deux positions $f$,$g$ de $!A$ si ces deux
  fonctions différent en un unique entier $i$ et si $(f(i),g(i)) \in
  E_A$ est un coup de $A$.\\
  Dans ce cas, on pose $\lambda_{!A}(f,g) = \lambda_A(f(i),g(i))$
\end{itemize}

Maintenant, essayons de calculer le comonoïde commutatif libre sur n'importe quel
objet $A \in \p$. On sait qu'à priori, $A$ n'a aucune chance d'être
pointé dans $\p$. Il faut donc travailler avec $Neg(A)$, qui lui est
pointé.

On peut alors définir $!Neg(A)$ comme ci-dessus. Nous avons besoin
maintenant d'un lemme auxiliaire pour montrer que cet objet est le bon
candidat au comonoïde.

\begin{lemma}
  Soit $(M,d,e)$ un comonoïde commutatif de $\p$.
  Alors $M$ est un jeu de Conway négatif
\end{lemma}

\begin{proof}
  Il suffit de regarder le diagramme
  $$
  \xymatrix{
    M \ar[r]^{d} \ar[rd]_{d} & M \tensor M \ar[d]^{\alpha_{M,M}}\\
    & M \tensor M 
  }
  $$
  qui indique que si $M$ possède un coup positif, alors la stratégie
  doit réagir en jouant ou bien à gauche, ou bien à droite. Or le
  morphisme $\alpha_{M,M}$ va renverser ce choix, ce qui empêche le
  diagramme précédant de commuter.
\end{proof}

\para{Remarque}
Ce lemme implique que tout comonoïde commutatif à une flèche unique
dans $I$ et donc les morphismes entre objets pointés coïncident avec
les morphismes, ce qui fait qu'on à bien construit le comonoïde
commutatif libre général.
\newline
 
Reste à démontrer la
proposition suivante:

\begin{prop}
  Soit $A \in \p$ quelconque.
  Alors $!Neg(A)$ est le comonoïde libre sur $A$
\end{prop}

\begin{proof}
  Soit un comonoïde commutatif $M$ tel qu'il existe un morphisme $f:M
  \rightarrow A$. 
  Alors, par l'adjonction entre $U$ et $Neg$, on obtient un morphisme
  de $Neg(f):M \rightarrow Neg(A)$. On en déduit qu'il existe un unique
  morphisme de comonoïde $!f$ qui fait commuter 
    $$
  \xymatrix{
    M \ar[r]^{Neg(f)} \ar[d]_{!Neg(f)} & Neg(A) \\
    {!Neg(A)} \ar[ur]_{d_{!Neg(A)}} & 
  }
  $$
  Mais l'adjonction nous permet de compléter le diagramme par 
  $$
  \xymatrix{
    M \ar[r]^{Neg(f)} \ar@/^1cm/[rr]^f \ar[d]_{!Neg(f)} & Neg(A) \ar[r]^{\mu} & A\\
    {!Neg(A)} \ar[ur] \ar[rru]_{d_{!A}}& &
  }
  $$
  On obtient donc l'existence et l'unicité du morphisme de comonoïde
  entre $M$ et $!Neg(A)$ qui fait commuter le diagramme ci-dessus.
\end{proof}

\para{Monoïde commutatif libre}
Comme le tenseur de la catégorie $\p$ possède aussi un adjoint à
droite, on peut de manière duale construire le monoïde commutatif
libre sur un objet $A$ en s'appuyant sur la catégorie des jeux
positifs. On obtient alors la construction duale de l'exponentielle, à
savoir $?A$.

\section*{Vers une théorie monoïdale}

Nous aimerions maintenant recomprendre ces travaux en termes
d'extensions de Kan, comme nous l'avons fait dans le cas où le tenseur
commute avec la somme. Il semble alors naturel de changer la catégorie
de base en travaillant avec $Inj$. Malheureusement, pour le moment,
nous n'arrivons pas à utiliser la monoïdalité de l'extension de Kan
par rapport à la convolution pour en déduire que cette extension
appliquée en $1$ nous donne le monoïde libre. Il faut comprendre
comment utiliser la commutation aux colimites filtrées pour pouvoir
exprimer la convolution en terme de tenseur comme dans la proposition
\ref{prop:convolution}. Il semble néanmoins possible d'établir des
encadrements de ces convolutions qui permettraient d'obtenir des
propriétés à la limite (car toute sous-catégorie infinie de la
catégorie $Inj$ voit son foncteur d'inclusion être final).

Cette dernière option semble néanmoins ouvrir la voie au développement d'une
théorie monoïdale à l'instar des théories algébriques de Lawvere. En
effet, malgré nos investigations, aucune théorie générale telle que les
catégories enrichies sur une catégorie monoïdale~\cite{kelly82} ne
paraît en mesure de capturer la notion de foncteur monoïdal nécessaire au
développement de la théorie susmentionnée. Il faut néanmoins investir
nos prochaines recherches vers la théorie des opérades qui semble
répondre partiellement à notre problématique.

\chapter{Un modèle de logique linéaire intuitionniste} \label{sct:asynchronous game}

Pour donner lieu à un modèle de logique linéaire intuitionniste, il
faut restreindre notre modèle à la partie négative des jeux de
Conway. Mais le problème qui se pose alors est de trouver une nouvelle
clôture car la catégorie n'est plus compacte close (le dual n'existe
plus). Heureusement, comme nous l'avons mentionné dans le chapitre
\ref{ch:intro}, un résultat nous permet d'exporter
automatiquement la fermeture via un adjonction monoïdale.

\section*{Petite construction catégorique}

Nous revenons maintenant sur une propriété énoncée au chapitre \ref{ch:intro}
\newline

\para{Une adjonction monoïdale donnant lieu à une clôture}
Soit $(\C,\tensor_\C,\pop_\C)$ une catégorie symétrique monoïdale fermée
et $(\D,\tensor_\D)$ un catégorie symétrique monoïdale.
Supposons qu'il existe un foncteur monoïdal fort $U:\D \rightarrow \C$
qui est à la fois plein et fidèle et qui possède un adjoint à droite
$F:\C \rightarrow \D$ formant l'adjonction $U \dashv F$.
    
Alors, on peut exporter la fermeture sur $\C$ en une fermeture sur
$\D$ en définissant pour $A,B$ dans $\D$:
$$
A \pop_\D B = F(U(A) \pop_\C U(B))
$$

Notons que la plupart du temps, le foncteur $U$ sera le foncteur
d'inclusion et $F$ forcera le fermeture de $\C$ à vivre dans la
sous-catégorie $\D$.

\begin{proof}
$$
\begin{array}{rcll}
  \D(B,A \pop_\D C) & \cong & \D(B,F(U(A) \pop_\C U(C))) &\\
  & \cong & \C(U(B),U(A) \pop_\C U(C)) & 
  \quad \textrm{adjonction } U \dashv F\\ 
  & \cong & \C(U(A) \tensor_\C U(B),U(C)) & 
  \quad \textrm{adjonction } \tensor_\C \dashv \pop_\C\\ 
  & \cong & \C(U(A \tensor_\D B),U(C)) & \quad
  \textrm{monoïdalite forte}\\
  & \cong & \D(A \tensor_\D B,C) & \quad \textrm{plein et fidèle}
\end{array}
$$
\end{proof}

\para{Un exemple en logique intuitionniste}
Soit $\C$ la catégorie des formules de logique linéaire, et $\D$
le fragment positif (ie. toutes les formules qui sont des tenseurs
d'exponentielles). Il est connu que $!$ est adjoint à droite de
l'inclusion et donc la clôture de $\C$ peut être exportée en:
\begin{equation*}
A \Rightarrow B \ = \ ! \ (A \pop B)
\end{equation*}

Cela donne lieu à l'un des deux codages de la logique intuitionniste,
l'autre étant le fragment négatif avec le flèche de co-Kleisli $!A
\pop B$. 

\section*{Jeux de Conway négatifs}

Un jeu de Conway négatif (à gain) est simplement un jeu de Conway (à
gain) dont tous les coups initiaux sont négatifs. On note $\n$ la
sous-catégorie des jeux négatifs. Le produit tensoriel peut être
directement exporté car il préserve la polarité des coups
initiaux. Nous allons utiliser le résultat ci-dessus pour obtenir une
clôture car la clôture de $\p$ n'est évidemment pas préservée, notre
catégorie négative n'ayant pas de dual.  

Pour cela, observons que l'inclusion, qui est toujours un foncteur
monoïdal fort plein et fidèle, admet ici un adjoint à droite en la
présence du foncteur $Neg : \p \rightarrow \n$ qui oublie les parties
commençant par des coups positifs.  
En effet, il est aisé de vérifier que 
$$
\frac{U A \rightarrow B}{A \rightarrow Neg(B)}
$$
définit bien une adjonction. On a alors directement la clôture donnée
par 
$$A \pop B = Neg(A^* \tensor B)$$ 

Nous pouvons maintenant nous occuper du produit cartésien.
\newline

\para{Produit}
Le produit de deux jeux négatifs $A$ et $B$, noté $A \&
B$ est défini comme:

\begin{itemize}
\item
  l'ensemble de ses positions est l'union disjointe des positions de
  $A$ et $B$ dans laquelle on a identifié les deux racines $\star_A$
  et $\star_B$ en la nouvelle racine $\star_{A \& B}$ de $A \& B$.
\item
  les coups (nécessairement opposants) partant de la racine sont de
  deux sortes:
  $$
  \star_{A \& B} \rightarrow 
  \left\{
  \begin{array}{l}
    x \mbox{ if } (\star_A,x) \in E_A \\
    y \mbox{ if } (\star_B,y) \in E_B \\
  \end{array}
  \right.
  $$
\item
  les coups partant d'une position $x$ de la composante $A$
  (resp. $B$) sont exactement les coups partant de $x$ dans $A$
  (resp. $B$) avec la même polarité
\item
  le gain d'un chemin $s$ dans la composante $A$ (resp. $B$) est le
  gain de ce chemin dans $A$ (resp. $B$)
\end{itemize}

Il est aisé de vérifier que cette définition satisfait aux propriétés
du produit.
\newline

\para{exponentielle}
Nous avons donné au chapitre \ref{ch:exponentielle} une construction
de l'exponentielle pour les jeux de Conway. Nous allons vérifier ici
que la même définition donne toujours le comonoïde commutatif libre
lorsqu'on se place dans la catégorie des jeux à gain.

\begin{prop}
  L'objet ${!A} = S^n(A)$ est le comonoïde commutatif libre sur $A$
  dans la catégorie de jeux de Conway à gain
\end{prop}

\begin{proof}

On doit regarder si toute les étapes de la construction donnent des
stratégies gagnantes.
 
Le seul point délicat à vérifier est que $!A$ est toujours la limite du
diagramme $\Delta_\omega \xrightarrow{T'} \p$, toutes les autres
constructions donnant trivialement lieu à des stratégies gagnantes. 

Il faut donc montrer pour tout cône sur $X$ que l'unique flèche de $X
\rightarrow {!A}$ faisant commuter les cônes est une stratégie
gagnante. 
  
Procédons par l'absurde et supposons que cette stratégie joue un
chemin $t$ dans une partie $s$ avec $\kappa^+(t) = 0 \wedge
\kappa^-(t) > 0$. Comme $s$ est nécessairement fini, il existe $n$ tel
que pour tout $n'>n$, toute les positions de $s$ (qui sont des
fonctions) envoient $n'$ sur $\star_A$. En d'autres termes, on n'a
joué que sur le jeu $X \rightarrow A^{\tensor n}$.
Or, le diagramme suivant commute
$$
\xymatrix{
{!A} \ar[r] \ar[rd] & X \ar[d]\\
& A^{\tensor n} 
}
$$
On en déduit $X \rightarrow A^{\tensor n}$ est perdante, et donc que
le cône de départ n'était pas pris sur les jeux de Conway à gain.
\end{proof}

\section*{Catégorie linéaire et modèle de LLI}

Nous voulons maintenant montrer que notre catégorie définit un modèle
correct de LLI. Pour cela, nous allons passer par la notion de
catégorie de Lafont

\begin{definition}[Catégorie de Lafont]
  Une catégorie de Lafont consiste en 
  \begin{itemize}
  \item
    une catégorie symétrique monoïdale close avec produits finis
    $(\C,\tensor,1,\&,\top)$ 
  \item
    pour tout objet $A\in \C$, l'objet $!A$ est le comonoïde
    commutatif libre sur $A$ 
  \end{itemize}
\end{definition}

Les catégories de Lafont possèdent malgré leur simplicité la propriété
que nous recherchons.

\begin{prop}
  Toute catégorie de Lafont induit un modèle correct de logique
  linéaire intuitionniste. 
\end{prop}

Or, tout marche bien dans notre cadre. En effet, la catégorie $\n$ est
symétrique monoïdale close avec produits finis. Elle possède de plus
un comonoïde commutatif libre sur chacun de ses objets. On en déduit
la proposition suivante 

\begin{prop}
  La catégorie des jeux de Conway négatifs à gain $\n$ est une
  \emph{catégorie de Lafont}.
\end{prop}

On en déduit alors le théorème qui a motivé tout ce chapitre

\begin{theoreme}
   La catégorie des jeux de Conway négatifs à gain $\n$ fournit un
   modèle correct de la logique linéaire intuitionniste.
\end{theoreme}
\para{Remarque}

La catégorie $\n$ ne nécessite pas de notion de gain pour être une
catégorie de Lafont. Néanmoins, c'est l'opportunité offerte par cette
notion dans la modélisation du contrôle qui a motivé une bonne
partie de notre travail. Il ne faut donc pas perdre à l'esprit que
même si le gain n'apparaît pas dans le prochain chapitre, il reste
fondamental dans notre approche de la sémantique des langages de
programmation. 
\newline

Nous avons maintenant tous les outils prérequis à la construction du
modèle de références globales et locales qui nous avons en tête.
Nous allons laisser de côté la notion de gain car le contrôle
n'apparaît pas encore dans notre langage.  

\chapter{Extension à un modèle de référence utilisant la trace}

Comme nous nous sommes concentrés sur la catégorie des jeux négatifs, il
est naturel de concevoir un langage en appel par nom. En effet, une
partie jouée par une stratégie entre deux jeux négatifs commence
toujours à droite; on peut donc composer une stratégie perdante
$\sigma:1 \rightarrow A$ avec une stratégie gagnante $\tau:A \rightarrow B$
et obtenir une stratégie gagnante $\sigma;\tau$. Ceci vient du fait
que $\tau$ peut se désintéresser de son argument, et ainsi lui donner
un mauvais argument en entrée ne la dérange pas. C'est typiquement le
genre d'interaction qui a lieu dans un langage en appel par nom.

Il est à noter que nous pourrions aussi bien décrire un langage type PCF
en appel par valeur en nous tournant vers la catégorie des jeux
positifs. À terme, nous voulons d'ailleurs décrire un cadre où appel par nom et
appel par valeur vivent dans la même catégorie et où passer de l'un à
l'autre est transparent.
Il nous semble en effet que la distinction appel par valeur/nom n'est
pas fondamentale et que les mêmes outils algébriques peuvent expliquer
ces deux cadres de manière naturelle.  

\section*{TracedAlgol: un langage avec références globales et locales} 

Comme nous l'avons déjà indiqué, nous donnons ici un langage avec
références mais sans aliasing. Nous devons donc avoir deux notions de
types, les types valeurs $A,B$ et les types références $\alpha,\beta$,
pour pouvoir distinguer une valeur et référence grâce au typage.
$$
A,B ::= \Unit \sep \Bool \sep \Nat \sep A \rightarrow B 
\sep A \times B \qquad \qquad \alpha,\beta ::= A \sep \Ref{A}
$$ 

\para{Remarque}
Cette restriction des références (qui ne peuvent pas pointer sur
d'autres références) permet tout de même de typer des programmes
d'ordre supérieur arbitraire. En revanche, on ne rentre pas dans le
cadre de l'aliasing, ce qui ferait entrer le langage dans une classe
de complexité bien plus grande.
\newline

Les termes du langage sont inspirés d'Algol:
\begin{eqnarray*}
M,N & ::= & \Skip \sep b \sep n \sep x \sep \lambda x.M \sep MN
\sep x := M \sep !x \sep \soit x:=M \dans N \\
& & | \zero(M) | \si M \alors M_1 \sinon M_2 \sep \langle M,N \rangle \sep \pi_1(M)
\sep \pi_2(M) \sep M;N 
\end{eqnarray*}
où $b$ vaut pour les booléens $T$ et $F$, $n$ représente un entier et 
$x$ est une variable.
\newline

\para{Remarque: $\mathrm{new}$ ou $\mathrm{let}$}
Nous notons ici le traditionnel terme $\mathrm{let} \ x:=M \dans N$ par
$\mathrm{new} \ x:=M \dans N$ pour rappeler que cette interprétation des
variables locales nous vient du travail de Milner sur la sémantique du
$\nu$ \cite{milner:ac5}.
\newline

Nous divisons le contexte d'interprétation en un contexte de variables
$\Gamma$ qui associe à chaque variable un type et un contexte de
référence $\Delta$ qui associe à chaque référence un type référence.
Nous supposerons toujours que $\mathrm{dom}(\Gamma) \cap
\mathrm{dom}(\Delta) = \emptyset$. 

Les règles de typage sont données en figure~\ref{fig:typage}

\begin{figure}[h]
  \begin{center}
    $$
    \rules{Var}{}{\Gamma,x:A; \  \vdash x:A}
    \qquad
    \rules{Bool}{}{\Gamma ; \ \vdash b : \Bool}
    \qquad
    \rules{Nat}{}{\Gamma ; \ \vdash n : \Nat}
    \qquad
    \rules{Unit}{}{\Gamma ; \ \vdash \Skip : \Unit}
    $$
    \vspace{0.2cm}
    $$
    \rules{Abs}{\Gamma,x:\alpha;\Delta \vdash M:\beta}{\Gamma;\Delta
    \vdash \lambda x.M : \alpha \rightarrow \beta} \, (x \notin \Gamma)
    \qquad
    \rules{App}{\Gamma;\Delta \vdash M:\alpha \rightarrow \beta \quad
    \Gamma;\Delta \vdash N : \alpha}{\Gamma;\Delta \vdash MN:\beta}
    $$
    \vspace{0.2cm}
    $$
    \rules{Trace}{\Gamma;\Delta \vdash M:A \quad \Gamma ; \Delta,
    x:\Ref{A} \vdash N:\beta}{\Gamma;\Delta \vdash \soit x:=M \dans N  :
    \beta} \, (x \notin \Delta)
    \qquad
    \rules{Weak}{\Gamma;\Delta \vdash
    M:\alpha}{\Gamma;\Delta , x:\Ref{A} \vdash M:\alpha}
    $$
    \vspace{0.2cm}
    $$
    \rules{Seq}{\Gamma;\Delta \vdash M:\Unit \quad \Gamma;\Delta
    \vdash N:\Unit}{\Gamma;\Delta \vdash M;N:\Unit}
    \qquad
    \rules{Zero}{\Gamma;\Delta \vdash M : \Nat}{\Gamma;\Delta \vdash
    \zero(M) : \Bool}
    $$
    \vspace{0.2cm}
    $$
    \rules{If}{\Gamma;\Delta \vdash M:\Bool \quad \Gamma;\Delta \vdash
    M_i : \alpha (i=1,2)}{\Gamma;\Delta \vdash \si M \alors M_1 \sinon
    M_2 : \alpha}
    $$
    \vspace{0.2cm}
    $$
    \rules{Pair}{\Gamma;\Delta \vdash
    M_i : \alpha_i (i=1,2)}{\Gamma;\Delta \vdash
    \langle M_1,M_2 \rangle : \alpha_1 \times \alpha_2}
    \qquad
    \rules{Proj}{\Gamma;\Delta \vdash
    M : \alpha_1 \times \alpha_2}{\Gamma;\Delta \vdash
    \pi_i(M) : \alpha_i (i=1,2)} 
    $$
    \vspace{0.2cm}
    $$
    \rules{Assign}{\Gamma;\Delta,x:\Ref{A} \vdash M:A}{\Gamma ;
    x:\Ref{A},\Delta \vdash x := M : \Unit} 
    \qquad
    \rules{Deref}{}{\Gamma ; x:\Ref{A},\Delta \vdash \ !x : A}
    $$
  \end{center}
  \caption{Règles de typage dans TracedAlgol}\label{fig:typage}
\end{figure}

Nous allons maintenant donner la sémantique opérationnelle de notre
langage. Celle-ci sera donnée pour des couples $(M,\sigma)$ appelés
configurations, où $M$ est un programme et $\sigma$ est une fonction
qui va des références dans les valeurs, que nous appellerons état de
la mémoire.
Nous faisons le choix de présenter une sémantique \emph{big step} car
elle forme un cadre plus agréable pour l'appel par nom.

Pour cela, il faut une notion de forme canonique, qui sont les termes
de la forme:
$$
V ::= n \sep b \sep x \sep \Skip \sep \lambda x.V \sep \langle V_1,V_2
\rangle
$$

Lorsque qu'une configuration $(M,\sigma)$ se réduit en une
configuration canonique $(V,\sigma')$, nous notons
$$
(M,\sigma) \reduct (V,\sigma')
$$ 

La figure \ref{fig:operation} décrit inductivement la sémantique
opérationnelle de notre langage.

\begin{figure}[h]
  \begin{center}
    $$
    \rules{Value}{}{(V,\sigma) \reduct (V,\sigma)}
    \qquad
    \rules{Lambda}{(M,\sigma) \reduct (\lambda x.V',\sigma') \quad
    (V'[N/x],\sigma') \reduct (V,\sigma'')}{(MN,\sigma) \reduct (V,\sigma'')}
    $$
    \vspace{0.2cm}
    $$
    \rules{Cond_T}{(M,\sigma) \reduct (T,\sigma') \quad (M_1,\sigma')
      \reduct (V,\sigma'')}{(\si M \alors M_1 \sinon M_2,\sigma) \reduct
      (V,\sigma'')}
    \qquad
    \rules{Cond_F}{(M,\sigma) \reduct (F,\sigma') \quad (M_2,\sigma')
      \reduct (V,\sigma'')}{(\si M \alors M_1 \sinon M_2,\sigma) \reduct
      (V,\sigma'')}
    $$
    \vspace{0.2cm}
    $$
    \rules{Seq}{(M,\sigma) \reduct (\Skip,\sigma') \quad (N,\sigma')
    \reduct (\Skip,\sigma'')}{(M;N,\sigma) \reduct (\Skip,\sigma'')}
    \qquad
    \rules{Zero_F}{(M,\sigma) \reduct (n+1,\sigma')}{(\zero(M),\sigma) \reduct
    (F,\sigma'')} 
    $$
    \vspace{0.2cm}
    $$
    \rules{Zero_T}{(M,\sigma) \reduct (0,\sigma')}{(\zero(M),\sigma) \reduct
    (T,\sigma'')} 
    \qquad
    \rules{Trace}{(M,\sigma) \reduct (V',\sigma') \quad
    (N,\sigma'\cup( x \mapsto V')) \reduct (V,\sigma'')}{(\soit x:=M
    \dans N,\sigma) \reduct (V,\sigma''\setminus x)}
    $$
    \vspace{0.2cm}
    $$
    \rules{Assign}{(M,\sigma) \reduct (V,\sigma')}{(x:=M,\sigma)
    \reduct (skip,\sigma' \cup (x \mapsto V)) }
    \qquad
    \rules{Deref}{\sigma(x) = V}{(!x,\sigma) \reduct (V,\sigma)}
    $$    
  \end{center}
  \caption{Sémantique opérationnelle de TracedAlgol}\label{fig:operation}
\end{figure}

\section*{Interprétation dans la catégorie des jeux de Conway négatifs}

Nous devons maintenant montrer comment nous interprétons ce langage
dans la catégorie de jeux de Conway négatifs à gain $\n$.

Chaque type valeur est interprété par un objet $\sem{A} \in \n$ comme
suit :
\begin{itemize}
  \item 
    $\sem{\Bool} = \bool ool$
  \item 
    $\sem{\Nat} = \mathbb{N}at$
  \item 
    $\sem{\Unit} = 1$
  \item
    $\sem{A \times B} = \sem{A} \& \sem{B}$
  \item
    $\sem{A \rightarrow B} = !\sem{A} \pop \sem{B} $
\end{itemize}

Le statut d'un type référence est un peu différent car il doit pouvoir
à la fois se situer dans le co-domaine et le domaine du programme, on
en déduit l'interprétation du jugement de type 
$$
\sem{x_1:A_1,\ldots,x_n:A_n;y_1:\Ref{B_1},\ldots,y_m:\Ref{B_m} \vdash
  M:\alpha}
$$ 
comme un stratégie 
$$
\sem{M} : !\sem{A_1} \tensor \ldots \tensor !\sem{A_n} \tensor
!\sem{B_1} \tensor \ldots \tensor !\sem{B_n} \rightarrow !\sem{B_1} \tensor
\ldots \tensor !\sem{B_n} \tensor \sem{M}
$$

Donnons maintenant l'interprétation des règles de typage. 
\begin{itemize}
\item
  Les règles pour les constantes sont évidentes
\item
  $Var$ correspond à l'identité (la stratégie \emph{copycat})
\item
  $Abs$ vient juste de la clôture de notre catégorie
  $$
  \sem{\Gamma;\Delta \vdash \lambda x.M : \alpha \rightarrow \beta} = \Lambda
  \sem{\Gamma,x:\alpha;\Delta \vdash M:\beta} : \sem{\Gamma}
  \rightarrow !\sem{A} \pop \sem{B}
  $$
\item
  $App$ vient de l'évaluation (la co-unité de la clôture) et de la co-multiplication
  de l'exponentielle
  $$
  \sem{\Gamma} \tensor \sem{\Delta} \xrightarrow{d}
  (\sem{\Gamma} \tensor \sem{\Delta}) \tensor (\sem{\Gamma} \tensor
  \sem{\Delta}) \rightarrow (\sem{\alpha} \pop
  \sem{\beta}) \tensor \sem{\alpha} \xrightarrow{eval} \sem{\beta}
  $$
\item
  $Seq$ : le terme $M;N$ correspond juste à du sucre syntaxique pour
  $\soit x:=M \dans N$ où $x$ n'apparaît pas dans $N$ 
\item
  $Weak$ correspond à l'identité car lorsqu'on ajoute une référence,
  elle apparaît des deux côtés du séquent
\item
  $Pair$ et $Proj$ sont décrites par les morphismes liés au produit
  cartésien
\item
  $If$ est donnée par la stratégie qui interroge son argument, exécute
  le premier programme si $V$ et le deuxième si $F$ 
\item
  $Zero$ est interprétée par la stratégie évidente
\item
  $Assign$ va être interprétée en deux temps :
  \begin{enumerate}
  \item
    d'abord effacer la valeur courante de la cellule,
  \item
    ensuite remplacer cette valeur par l'interprétation de $M$
  \end{enumerate}
  Plus précisément, le terme $M$ est interprété par un morphisme
  $$
  \sem{\Gamma} \tensor
  !\sem{A} \tensor \sem{\Delta} \rightarrow !\sem{A'} \tensor
  \sem{\Delta} \tensor \sem{A''} 
  $$
  (Notons que $A$,$A'$ et $A''$ coïncident dans le modèle, nous les
  distinguons uniquement pour faciliter la lecture de la construction)
 
  Comme tous les objets autres que $\sem{A''}$ sont des $!$-coalgèbre,
  on peut alors appliquer le foncteur $!$ et la co-multiplication de
  $!$ vue comme une comonade. On obtient 
  $$
  \sem{\Gamma} \tensor
  !\sem{A} \tensor \sem{\Delta} \rightarrow !\sem{A'} \tensor
  \sem{\Delta} \tensor !\sem{A''} 
  $$
  (on vient juste d'exprimer la règle de promotion en logique linéaire)

  On interprète $x:=M$ en post-composant avec la co-unité en $A'$ puis
  la permutation de $\Delta$ et $A''$
  $$
  !\sem{A'} \tensor \sem{\Delta} \tensor !\sem{A''} \rightarrow
  \sem{\Delta} \tensor !\sem{A''} \rightarrow !\sem{A''} \tensor
  \sem{\Delta} 
  $$  
  On obtient ainsi le morphisme 
  $$
  \sem{\Gamma} \tensor
  !\sem{A} \tensor \sem{\Delta} \rightarrow !\sem{A''} \tensor
  \sem{\Delta}   
  $$
  qui interprète le séquent 
  $$
  \Gamma ;
  x:\Ref{A},\Delta \vdash x := M : \Unit
  $$
\item
  $Deref$ correspond juste à la duplication de la référence par la
  co-multiplication puis une permutation (encore une fois, nous
  distinguons artificiellement les trois copies $A$ par souci de
  clarté)  
  $$
  \sem{\Gamma} \tensor
  !\sem{A} \tensor \sem{\Delta} \rightarrow !\sem{A'} \tensor
  \sem{\Delta} \xrightarrow{d_{!\sem{A'}}} !\sem{A'} \tensor !\sem{A''} \tensor
  \sem{\Delta} \xrightarrow{\alpha_{!\sem{A''},\sem{\delta}}} !\sem{A'} \tensor \sem{\Delta} \tensor
  !\sem{A''}
  $$
\item
  $Trace$ : c'est là que réside l'originalité de notre interprétation.
  Comme pour l'assignation, on construit, à partir de l'interprétation
  de $M$, le terme
  $$
  \sem{\Gamma} \tensor \sem{\Delta} \rightarrow \sem{\Delta} \tensor
  !\sem{A''} 
  $$
  que l'on ``tensorise'' avec l'identité $!\sem{A} \rightarrow
  !\sem{A'}$ (plus un peu de commutation)
  $$
  \sem{\Gamma} \tensor !\sem{A} \tensor \sem{\Delta} \rightarrow
  !\sem{A'} \tensor \sem{\Delta} \tensor !\sem{A''} 
  $$
  À nouveau, on applique la même technique que pour l'assignation en
  post-composant avec la co-unité en $A'$ puis la permutation de
  $\Delta$ et $A''$, et on obtient ainsi une interprétation de la
  création de la référence $x$ stockant la valeur $M$
  $$
  \sem{\Gamma} \tensor !\sem{A} \tensor \sem{\Delta} \rightarrow
  !\sem{A''} \tensor \sem{\Delta}
  $$
  On compose ensuite avec l'interprétation de $N$, et on a
  $$
  \sem{\Gamma} \tensor
  !\sem{A} \tensor \sem{\Delta} \rightarrow !\sem{A''} \tensor
  \sem{\Delta} \tensor \sem{\beta}
  $$
  Il ne reste plus qu'à tracer sur $!\sem{A}$ (après permutation) et
  on obtient l'interprétation souhaitée du terme
  $$
  \soit x:=M \dans N : \sem{\Gamma} \tensor \sem{\Delta}
  \rightarrow \sem{\Delta} \tensor \sem{\beta}
  $$
\end{itemize}

\para{Remarque}
Il est évident que l'interprétation donnée est stable par contexte.

\section*{Correction Équationnelle}

Nous décrivons d'abord la notion usuelle d'équivalence observationnelle
pour laquelle nous voulons un résultat de correction du modèle.

Soit $M$ un terme clos et sans emplacement mémoire libre
de TracedAlgol. On note $M\reduct$ si $M$ est typable et $(M,\emptyset) \reduct
(V,\sigma)$ pour un certain terme $V$.

\begin{definition}[équivalence observationnelle]
  Soit $M$ et $N$ deux termes. Alors $M$ et $N$ sont équivalents
  observationnellement, noté $M \simeq N$, ssi pour tout contexte
  $C[-]$ tel que $C[M]$ et $C[N]$ sont clos et sans emplacement
  mémoire libre, on a $C[M] \reduct$ ssi $C[N] \reduct$ 
\end{definition} 

Il est notable que la condition $C[M] \reduct$ ssi $C[N] \reduct$
suffise à s'assurer que, si $M$ et $N$ se réduisent en des valeurs,
celles-ci sont identiques (grâce au test en zéro et au $if$). 

Pour alléger la notation dans ce qui suit, nous noterons $\soit
x:=V_1,y:=V_2 \dans M$ pour $\soit x:=V_1 \dans (\soit y:=V_2 \dans
M)$.
De la même manière, étant donné un état de la mémoire $\sigma$, on
note $\soit \sigma \dans M$ le terme $\soit x_1:=\sigma(x_1),
\ldots,x_n:=\sigma(x_n) \dans M$ où les $x_i$ parcourent les états
interprétés par $\sigma$.  

Dans la lignée des preuves de correction en sémantique des jeux
\cite{ahm}, nous décomposons la preuve de correction équationnelle en
deux étapes : correction et adéquation.

\begin{lemma}[correction]
  Soit $M$ un terme. Si $(M,\sigma) \reduct (V,\sigma')$, alors
  $\sem{\soit \sigma \dans M} = \sem{\soit \sigma' \dans V}$
\end{lemma}

\begin{proof}
  On opère par une induction standard sur la dérivation de $(M,\sigma)
  \reduct (V,\sigma')$ en utilisant les équations de la
  figure~\ref{fig:eqRef} et la relation
  $$
  \sem{\Gamma;\Delta \vdash \soit x:=V \dans (\lambda y.M)(!x) } =
  \sem{\Gamma;\Delta \vdash \soit x:=V \dans (\lambda y.M) V}
  $$
  dont la validité est assurée en étudiant précisément le comportement
  de la stratégie d'évaluation vis-à-vis des références.
\end{proof}

\begin{figure}[h]
  \begin{eqnarray*}
    \sem{\Gamma;\Delta \vdash \soit x:=V_1,y:=V_2 \dans M} & = &
    \sem{\Gamma;\Delta \vdash \soit y:=V_2,x:=V_1 \dans M}
    \\
    \sem{\Gamma;x:\Ref{A},\Delta \vdash \soit y:=V_2 \dans x:=V_1;M} & = &
    \sem{\Gamma;x:\Ref{A},\Delta \vdash x:=V_1;\soit y:=V_2 \dans M}
    \\
    \sem{\Gamma;\Delta \vdash \soit x:=V_1,x:=V_2 \dans M} & = &
    \sem{\Gamma;\Delta \vdash \soit x:=V_2 \dans M}  
  \end{eqnarray*}
  \caption{Équations concernant les références locales}\label{fig:eqRef}
\end{figure}

Nous avons besoin d'une sorte de réciproque appelé adéquation

\begin{lemma}[Adéquation]
Pour tout terme clos $M$, si $\sem{M} \neq \bot$ alors $M \reduct$
\end{lemma}

La preuve s'appuie sur une longue étude des parties d'une stratégie non
vide. De telles parties proviennent de l'interaction de plusieurs
parties provenant des sous termes de $M$. La taille de cette
interaction forme le limon d'une récurrence.

Il nous est maintenant possible de statuer sur la correction de notre
modèle

\begin{theoreme}[Correction Équationnelle]
  Si $M$ et $N$ sont des termes du même type et $\sem{M} = \sem{N}$,
  alors $M \simeq N$
\end{theoreme}

\begin{proof}
  Soit $C[-]$ tel que $C[M]$ et $C[N]$ sont clos et sans emplacement
  mémoire libre.
  Par symétrie, il nous suffit de montrer que $(C[M] \reduct)
  \Rightarrow (C[N] \reduct)$.
 
  $C[M] \reduct$ signifie qu'il existe $\sigma,V$ tel que
  $(C[M],\emptyset) \reduct (V,\sigma)$.
  On déduit de la correction que $\sem{C[M]} = \sem{\soit \sigma \dans
  V}$
  Comme l'interprétation est stable par contexte, on a nécessairement
  $\sem{C[N]} = \sem{\soit \sigma \dans V}$.

  Ceci donne une interprétation non vide pour $C[N]$, de laquelle on
  déduit via l'adéquation 
  $C[N] \reduct$
  
\end{proof}

\chapter*{Conclusion et travaux futurs}

Lors de ce stage de DEA (qui dura de manière inhabituelle un peu plus d'un
an), nous nous sommes fixé pour but la compréhension sémantique des
langages avec références à travers un cadre algébrique. Ceci dans l'espoir
d'étendre l'isomorphisme de Curry-Howard aux langages de programmation
impératifs.


Bien évidemment, ce programme ambitieux n'est pas encore abouti, mais
nous avons tout de même développé plusieurs outils importants nous
permettant de nous en approcher.
\begin{enumerate}
\item
  Nous avons construit un modèle de sémantique des jeux parenthésés de
  logique linéaire intuitionniste disposant de plus d'un opérateur de
  trace. Pour cela, nous avons utilisé le modèle des jeux de Conway,
  augmenté avec une notion de gain définie de manière axiomatique,
  rapprochant le gain à la notion de distance dans un espace
  géométrique. Ce modèle très riche peut être à la base de nombreux
  travaux sémantiques car les intuitions usuelles en matière
  d'opérateurs logiques ou catégoriques s'y expriment fort bien.
\item
  Nous avons donné un cadre catégorique pour la construction du
  comonoïde commutatif libre, que nous avons utilisé pour obtenir
  l'exponentielle sur les jeux de Conway à gain. Cette technique
  devrait ouvrir la porte à de multiples constructions de
  l'exponentielle dans des cadres où sa définition ``à la main'' n'est
  pas toujours accessibles. On pourra aussi s'en servir pour mieux
  comprendre son existence et sa définition dans certains cadres
  sémantiques. 
\item
  Nous avons utilisé notre nouveau cadre pour décrire un modèle
  d'un langage de type Algol avec fonctionnelle d'ordre supérieur.
  Le modèle défini dépend plus des propriétés catégoriques qu'il
  vérifie, que de ses propriétés intrinsèques. Nous pouvons donc
  considérer que nous avons décrit un modèle catégorique de notre
  langage avec référence.
\end{enumerate}

\vspace{1cm}

Dans un futur proche, nous voulons recomprendre notre construction de
l'exponentielle en terme d'extension de Kan dans une théorie
monoïdale. 
Nous avons aussi pour projet de fusionner notre modèle avec celui des
jeux asynchrones afin d'avoir en notre possession toute la puissance
de ce cadre.
Enfin, nous espérons utiliser notre travail pour obtenir un modèle de
langage impératif avec aliasing, et même avoir un cadre
sémantique unifié capturant aussi bien langages de haut que de bas
niveau.


\bibliographystyle{apalike}
\bibliography{rapport_dea.bib}

\def\figurename{\textsc{Fig. A.} \hspace{-0.265cm}}

\chapter*{Annexes}

Les lecteurs n'étant pas familié avec la théorie des catégories sont
rapportés à l'incontournable référence \cite{macLane71} pour la
découverte de ce magnifique champ mathématique. Nous présentons ici
les notions les moins usuelles abordées dans cette recherche. 

\section*{Catégorie monoïdale tracée}\label{section/tmc}

Une catégorie symétrique monoïdale tracée
\cite{joyal-street-verity} est une SMC $(\cat,\tensor,I,s)$ munie
d'une famille de fonction   
$$
Tr_X: \qquad \frac{X \tensor A \quad \longrightarrow \quad X \tensor B}
{A \quad \longrightarrow \quad B}
$$

vérifiant les axiomes suivants :

\begin{itemize}
\item {\bf Naturality :}
    $Tr_X(id_X \tensor g ; f ; id_X \tensor h) = g ; Tr_X f ; h$ 
\item {\bf Strength :}
  $Tr_X(f \tensor g) = Tr_X f \tensor g$
\item {\bf Symmetry sliding :}
    $Tr_X(Tr_Y(f;c_{XY} \tensor id_B)) = Tr_Y(Tr_X(c_{XY} \tensor id_A ; f))$
\item {\bf Yanking :}
  $Tr_X (c_{XX}) = 1_X$
\end{itemize} 

\section*{Un exemple de catégorie monoïdale tracée}\label{sct:relation}

Nous mentionnons un exemple important de catégorie symétrique
monoïdale tracée. En effet, cette catégorie donne un opérateur de
trace pour la sémantique statique associée à la dynamique des jeux de
Conway. 

Soit la catégorie ${\bf Rel}$ des ensembles avec relations, 
munie du produit tensoriel défini sur les objets comme le produit
cartésien des ensembles et sur les relations par $\langle x,y
\rangle(R \times R') \langle x',y' \rangle$ ssi $x R y$ et $x' R y'$. 
Remarquons que ceci ne définit pas un produit au sens catégorique.
Pour $R : X \times A \rightarrow X \times B$, on définit $Tr_X
R : A \rightarrow B$ via :
$$
a (Tr_X R) b \equi \exists x \in X . \ \langle x,a \rangle R \langle x ,
b \rangle 
$$

On en déduit que $\langle {\bf Rel}, \times , Tr \rangle$ est une
catégorie tracée.

\section*{Extension de Kan}

Si $\A$ est un sous-ensemble de $\B$, une fonction $A:\A
\rightarrow\C$ vers un ensemble non vide $\C$ peut être étendue sur
$\B$ de beaucoup de manières, mais il n'y a pas de façon canonique de
le faire. Cependant, si $\A$ est une sous catégorie de $\B$, tout
foncteur $A:\A \rightarrow \C$ possède en principe deux extensions
canoniques (ou extrêmes) de $\A$ vers des foncteurs $L,R:\B
\rightarrow \C$. Ces extensions sont caractérisées par l'universalité
de transformations naturelles appropriées; elles n'existent pas
toujours mais on peut les calculer lorsque la catégorie $\A$ est
``petite'' et lorsque $\C$ est bi-complète. Nous ne présentons ici que
la définition de l'extension de Kan à gauche.

On part d'un foncteur $J:\A \rightarrow \B$ qui fait intuitivement de
la catégorie $\A$ une catégorie ``incluse'' dans $\B$.
Le problème est alors étant donné un foncteur $K:\A \rightarrow \C$ de
trouver l'extension naturelle $\exists_J A : \B \rightarrow \C$ munie
d'une transformation naturelle $\epsilon : A \rightarrow (\exists_J A)
\circ J $ universelle au sens où pour tout autre couple $S,\alpha : A
\rightarrow S \circ J$, il existe une unique transformation naturelle
$\sigma : \exists_J A \rightarrow S$ tel que $\alpha = \sigma J \cdot
\epsilon$

\begin{figure}[h!]
\begin{center}
\input{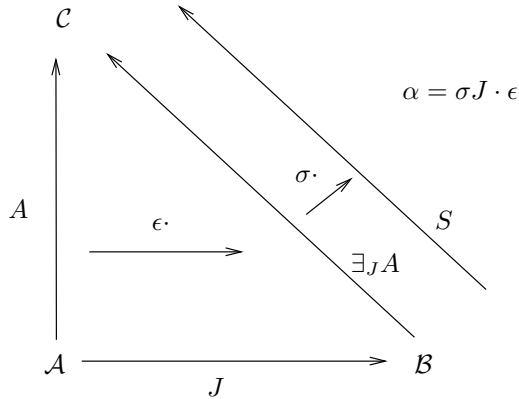}
\end{center}
\caption{Extension de Kan à gauche} 
\end{figure}

\section*{Convolution entre deux foncteurs d'une catégorie
  monoïdale dans un catégorie cocomplète}

Si $\C$ est une catégorie monoïdale cocomplète et $\A$ est un
catégorie monoïdale ``petite'', alors on peut munir la catégorie des
foncteurs $[\A,\C]$ d'une convolution monoïdale donnée par 
$$
F * G = \int^{A,A'} \A(A \tensor A',-) \tensor (FA \tensor GA')
$$

\section*{Catégories $\omega$-filtrées et colimites}

Une catégorie $\C$ est $\omega$-filtrée si 
\begin{itemize}
\item
  Pour tous objets $A,B$, il existe $C$ tel que 
  $$
  \xymatrix@=1pt{
    A \ar@{.>}[rrrrrrd] & & & & & & \\
    & & & & & & C \\
    B \ar@{.>}[urrrrrr] & & & & & & 
  }
  $$
\item
  Pour toutes flèches $f,g:A \rightarrow B$, il existe une flèche $h:B
  \rightarrow C$ tel que le diagramme
  $$
  \xymatrix@=1pt{
    & & & & & B \ar@{.>}[rrrrrd]^{h} & & & & & \\
    A \ar[rrrrru]^f \ar[rrrrrd]_g & & & & & & & & & & C \\
    & & & & & B \ar@{.>}[rrrrru]_{h} & & & & & 
  }
  $$
  commute.
\end{itemize}

\vspace{1cm}

Une colimite est dite \emph{filtrée} si elle est calculée sur une
catégorie $\omega$-filtrée. 

\vspace{1cm}

Traditionnellement, les colimites n'étaient calculées que sur de
préordres dirigés, qui ont ensuite été étendus à la notion de
catégories filtrées. Cette restriction s'est avérée finalement
inutile, mais le concept de colimites filtrées à garder son 
intérêt par la formule d'inversion de l'ordre d'application entre
colimites filtrées et limites finies. Dans notre travail, nous avons
regardé des catégories monoïdales dont le tenseur commute aux
colimites filtrées. 

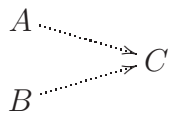
\end{document}